\documentclass[leqno,12pt,twoside]{amsart}

\usepackage{times}

\usepackage{amsmath, amssymb, amsfonts, latexsym, mdwlist}
\usepackage{subfig}
\usepackage{graphicx}
\usepackage{wrapfig}

\usepackage{paralist}
\setdefaultenum{(a)}{(i)}{}{}

\def\C{\ensuremath{\mathbb{C}}}

\def\H{\ensuremath{\mathbb{H}}}

\def\N{\ensuremath{\mathbb{N}}}
\def\P{\ensuremath{\mathbb{P}}}
\def\Q{\ensuremath{\mathbb{Q}}}
\def\R{\ensuremath{\mathbb{R}}}
\def\Z{\ensuremath{\mathbb{Z}}}

\newcommand{\dotcup}{\ensuremath{\mathaccent\cdot\cup}}

\def\Aut{\mathop{\mathrm{Aut}}\nolimits}

\def\ch{\mathop{\mathrm{ch}}\nolimits}

\def\Qcoh{\mathop{\mathrm{QCoh}}\nolimits}
\def\Coh{\mathop{\mathrm{Coh}}\nolimits}

\def\cok{\mathop{\mathrm{coker}}}

\def\dim{\mathop{\mathrm{dim}}\nolimits}

\def\inf{\mathop{\mathrm{inf}}\nolimits}
\def\End{\mathop{\mathrm{End}}}

\def\Ext{\mathop{\mathrm{Ext}}\nolimits}

\def\GL{\mathop{\mathrm{GL}}}
\def\Hom{\mathop{\mathrm{Hom}}\nolimits}

\def\RHom{\mathop{\mathbf{R}\mathrm{Hom}}\nolimits}
\def\id{\mathop{\mathrm{id}}\nolimits}
\def\Id{\mathop{\mathrm{Id}}\nolimits}
\def\im{\mathop{\mathrm{im}}\nolimits}

\def\Ker{\mathop{\mathrm{Ker}}\nolimits}

\def\min{\mathop{\mathrm{min}}\nolimits}

\def\Spec{\mathop{\mathrm{Spec}}}
\def\SL{\mathop{\mathrm{SL}}}
\def\ST{\mathop{\mathrm{ST}}\nolimits}

\def\Tot{\mathop{\mathrm{Tot}}\nolimits}

\newenvironment{Prf}{\textit{Proof.}\/}{\hfill$\Box$}

\def\MG13{\ensuremath{{\mathcal M}_{\Gamma_1(3)}}}
\def\tildeMG13{\ensuremath{\widetilde{\mathcal M}_{\Gamma_1(3)}}}
\def\Stab{\mathop{\mathrm{Stab}}}
\def\into{\ensuremath{\hookrightarrow}}
\def\onto{\ensuremath{\twoheadrightarrow}}

\def\blank{\underline{\hphantom{M}}}

\def\Db{\mathrm{D}^{\mathrm{b}}}

\def\pt{[\mathrm{pt}]}

\def\abs#1{\left\lvert#1\right\rvert}

\newcommand\stv[2]{\left\{#1\,\colon\,#2\right\}}

\newtheorem{Thm-s}[subsection]{Theorem}
\newtheorem{Prop-s}[subsection]{Proposition}
\newtheorem{Lem-s}[subsection]{Lemma}
\newtheorem{Thm}[subsection]{Theorem}
\newtheorem{Prop}[subsection]{Proposition}

\newtheorem{Lem}[subsection]{Lemma}
\newtheorem{Cor}[subsection]{Corollary}

\newtheorem{thm-int}{Theorem}

\theoremstyle{definition}

\newtheorem{Def}[subsection]{Definition}
\newtheorem{Rem}[subsection]{Remark}

\def\C{\ensuremath{\mathbb{C}}}

\def\H{\ensuremath{\mathbb{H}}}
\def\N{\ensuremath{\mathbb{N}}}
\def\P{\ensuremath{\mathbb{P}}}
\def\Q{\ensuremath{\mathbb{Q}}}
\def\R{\ensuremath{\mathbb{R}}}
\def\Z{\ensuremath{\mathbb{Z}}}

\def\AA{\ensuremath{\mathcal A}}

\def\CC{\ensuremath{\mathcal C}}
\def\DD{\ensuremath{\mathcal D}}
\def\EE{\ensuremath{\mathcal E}}
\def\FF{\ensuremath{\mathcal F}}
\def\GG{\ensuremath{\mathcal G}}
\def\HH{\ensuremath{\mathcal H}}

\def\KK{\ensuremath{\mathcal K}}

\def\NN{\ensuremath{\mathcal N}}
\def\OO{\ensuremath{\mathcal O}}
\def\PP{\ensuremath{\mathcal P}}

\def\QQ{\ensuremath{\mathcal Q}}
\def\TT{\ensuremath{\mathcal T}}
\def\UU{\ensuremath{\mathcal U}}
\def\WW{\ensuremath{\mathcal W}}

\def\YY{\ensuremath{\mathcal Y}}
\def\ZZ{\ensuremath{\mathcal Z}}

\def\AAA{\mathfrak A}

\def\DDD{\mathfrak D}
\def\EEE{\mathfrak E}
\def\FFF{\mathfrak F}

\def\MMM{\mathfrak M}
\def\RRR{\mathfrak R}
\def\SSS{\mathfrak S}

\newcommand{\ignore}[1]{}

\begin{document}
\title{The space of stability conditions on the local projective plane}

\author{Arend Bayer}
\address{Department of Mathematics, University of Connecticut U-3009, 196 Auditorium Road,
Storrs, CT 06269-3009, USA}
\email{bayer@math.uconn.edu}
\urladdr{http://www.math.uconn.edu/~bayer/}

\author{Emanuele Macr\`i}
\address{Department of Mathematics, University of Utah, 155 South 1400 East, Salt Lake City, UT 84112-0090, USA \& Mathematical Institute, University of Bonn, Endenicher Allee 60, D-53115 Bonn, Germany}
\email{macri@math.uni-bonn.de}
\urladdr{http://www.math.uni-bonn.de/~macri/}


\date{\today}

\begin{abstract}
We study the space of stability conditions on the total space of the canonical bundle over the projective plane. We
explicitly describe a chamber of geometric stability conditions, and show that its translates via autoequivalences cover
a whole connected component.  We prove that this connected component is simply-connected. We determine the group of
autoequivalences preserving this connected component, which turns out to be closely related to
$\Gamma_1(3)$.

Finally, we show that there is a submanifold isomorphic to the universal covering
of a moduli space of elliptic curves with $\Gamma_1(3)$-level structure. The morphism is
$\Gamma_1(3)$-equivariant, and is given by solutions of Picard-Fuchs equations. This result is motivated
by the notion of $\Pi$-stability and by mirror symmetry.
\end{abstract}

\maketitle

\section{Introduction}

In this paper, we study the space of stability conditions on the derived category of the local $\P^2$. Our approach is based on the chamber decomposition given by the wall-crossing for stable objects of the class of skyscraper sheaves of points.

\subsection{Motivation}
Consider a projective Calabi-Yau threefold $Y$ containing a projective plane $\P^2 \subset Y$. Ideally, one
would like to study the space of Bridgeland stability conditions on its derived category
$\Db(Y)$. Understanding the geometry of this space would give insights on the group of
autoequivalences of $\Db(Y)$ and give a global picture of mirror symmetry.  Understanding
wall-crossing for counting invariants of semistable objects would have many implication for
Donaldson-Thomas type invariants on $Y$.

However, no single example of stability condition on a projective Calabi-Yau threefold has been constructed. Instead,
in this article we focus on the full subcategory $\Db_{\P^2}(Y)$ of complexes concentrated on $\P^2$. The local model
for this situation is the total space $X = \Tot \OO_{\P^2}(-3)$ of the canonical bundle of $\P^2$, called the ``local
$\P^2$'': $\Db_{\P^2}(Y)$ is then equivalent to the derived category $\DD_0 := \Db_0(X)$ of coherent sheaves supported
on the zero-section.

Denote by $\Stab(\DD_0)$ the space of stability conditions $\sigma=(Z, \PP)$ on $\DD_0$ (see Appendix
\ref{app:BridgelandFramework} for a quick introduction to stability conditions). It is a three-dimensional
complex manifold coming with a local homeomorphism $\ZZ \colon \Stab(\DD_0) \to \Hom(K(\DD_0), \C) \cong \C^3$, $\ZZ((Z,\PP))=Z$.  The
goal of this article is to study the space $\Stab(\DD_0)$ as a test case for the properties we would expect in the case
of $Y$.

This space was first studied in \cite{Bridgeland:stab-CY}, where it was suggested that
the space is closely related to the Frobenius manifold of the quantum cohomology of $\P^2$. 
Further, understanding how Donaldson-Thomas type counting invariants of semistable objects 
depend on the stability conditions $\sigma \in \Stab(\DD_0)$ (i.e., wall-crossing phenomena) would be highly
interesting. For example, due to the derived equivalence $\DD_0 \cong \Db([\C^3/\Z_3])$ of \cite{Mukai-McKay} it
would give a new explanation for the relation between the Gromov-Witten potentials of $X$ and of $\C^3/\Z_3$
(``crepant resolution conjecture'', see \cite{Coates-CRC, CCIT:computing}). It could also explain the modularity
properties of the Gromov-Witten potential of $X$ observed in \cite{ABK:topological_strings}.

While these questions remain open, our results give a good description of a connected component
of $\Stab(\DD_0)$, explain its relation to autoequivalences of $\DD_0$, and do give a global mirror symmetry picture.

\subsection{Geometric stability conditions} \label{secti:geomstab}
In order to study $\Stab(\DD_0)$, we use one of its chamber decompositions.
We consider a chamber $U \subset \Stab(\DD_0)$ consisting of ``geometric'' stability conditions, which have the property
that all skyscraper sheaves $k(x)$, $x\in\P^2$, are stable of the same phase (see Definition \ref{def:GeomStability}
for the precise definition).

Our first result is a complete description of the geometric chamber (see Theorem
\ref{thm:geom-stability}): $U$ is an open, connected, simply-connected, 3-dimensional subset of $\Stab(\DD_0)$. Up to
shifts, a stability condition $(Z, \PP) \in U$ is determined by its central charge $Z$, and we give explicit
inequalities cutting out the set $\ZZ(U) \subset \Hom(K(\DD_0), \C) \cong \C^3$ of central charges $Z$
for $(Z, \PP) \in U$. The most interesting part of the boundary of $U$ has a fractal-like structure; its shape is
determined by the set of Chern classes of semistable vector bundles on $\P^2$.

Let $\Stab^\dag(\DD_0)$ be the connected component of $\Stab(\DD_0)$ containing $U$ and let $\overline{U}$ be the
closure of $U$ in $\Stab^\dag(\DD_0)$.  We can directly construct every \emph{wall} of $U$, i.e., the components
of  the boundary $\partial U=\overline{U}\setminus U$ of $U$ (see Theorem \ref{thm:boundary}).
We use this to prove the following result (see Corollary \ref{cor:ConnectedComponent}):

\begin{thm-int}\label{thmi:translates}
The translates of $\overline{U}$ under the group of autoequivalences generated by spherical twists at spherical sheaves in $\DD_0$ cover the whole connected component $\Stab^\dag(\DD_0)$.
\end{thm-int}
The translates of $U$ are disjoint, and each translate is a chamber on which the moduli space of stable
objects of class $[k(x)]$ is constant.

\subsection{Topology of $\Stab^\dag(\DD_0)$.}
In \cite{Bridgeland:stab-CY}, Bridgeland described an open connected subset $\Stab_a$ of $\Stab^\dag(\DD_0)$ consisting
of ``algebraic'' stability conditions that can be described in terms of quivers.
We will see that the subset $\Stab_a$ is not dense (in particular, it does not contain the ``large volume limit''
point: see Proposition \ref{prop:geomvsalg}). Nevertheless, by combining Bridgeland's description of
$\Stab_a$ with Theorem \ref{thmi:translates}, we prove the following result:

\begin{thm-int} \label{thmi:sc}
The connected component $\Stab^\dag(\DD_0)$ is simply-connected.
\end{thm-int}

\subsection{Autoequivalences}
In our situation, the local homeomorphism $\ZZ : \Stab^\dag(\DD_0) \to  \Hom(K(\DD_0), \C)$ is not a covering of its image.
This is a fundamental difference to the case of Calabi-Yau 2-categories (as studied in \cite{Bridgeland:K3,Thomas:stability,Bridgeland:ADE,Ishii-Ueda-Uehara, HMS:generic_K3s}).
Further, there is no non-trivial subgroup of autoequivalences of $\DD_0$ that acts as a 
group of deck transformation of the map $\ZZ$.  But, in any case, using Theorem 1 we can classify all autoequivalences
$\Aut^\dag(\DD_0)$
which preserve the connected component $\Stab^\dag(\DD_0)$:

\begin{thm-int} \label{thmi:auto}
The group $\Aut^\dag(\DD_0)$ is isomorphic to a product $\Z\times\Gamma_1(3)\times\Aut(\hat X)$.
\end{thm-int}

Recall that the congruence subgroup $\Gamma_1(3) \subset \SL(2, \Z)$ (see Section
\ref{sec:autoequivalences} for the definition) is a group on two generators $\alpha$ and $\beta$
subject to the relation $(\alpha\beta)^3 = 1$.  It is isomorphic to the subgroup generated by the
spherical twist at the structure sheaf $\OO_{\P^2}$ of the zero-section
$\P^2\into X$, and by the tensor product with $\OO_X(1)$.
The group $\Z$ is identified with the subgroup generated by the shift by $1$ functor $[1]$ and $\Aut(\hat X)$ denotes the group of automorphisms of the formal completion $\hat X$ of $X$ along $\P^2$.

\subsection{$\Pi$-stability and mirror symmetry}
Stability conditions on a derived category were originally introduced by Bridgeland in \cite{Bridgeland:Stab} to
give a mathematical foundation for the notion of $\Pi$-stability in string theory, in particular in Douglas' work, see e.g.
\cite{Douglas:stability, Aspinwall-Douglas:stability} and references therein. However, it has been understood that only
a subset of Bridgeland stability conditions is physically meaningful, i.e., there is a submanifold $M$ of the space of
stability conditions on $Y$ that parametrizes $\Pi$-stability conditions, and that is isomorphic to (the universal covering of) the complex K\"ahler moduli space. In fact, $M$ is (the universal covering of) a slice of the moduli space of SCFTs containing the sigma model associated to $Y$; in the physics literature, it is often
referred to as the ``Teichm\"uller space''.

By mirror symmetry, $M$ is also isomorphic to the universal covering of the moduli space of mirror partners
$\widehat{Y}$ of $Y$.  As explained in \cite{Bridgeland:spaces}, this leads to a purely algebro-geometric mirror
symmetry statement; we prove such a result in Section \ref{sec:MS}:

The mirror partner for the local $\P^2$ is the universal family over the moduli space $\MG13$ of elliptic curves
with $\Gamma_1(3)$-level structures. Its fundamental group is $\Gamma_1(3)$. Let
$\tildeMG13$ be the universal cover, with $\Gamma_1(3)$ acting as the group of deck transformations.

\begin{thm-int}\label{thmi:MS}
There is an embedding $I \colon \tildeMG13 \into \Stab^\dag(\DD_0)$ which is equivariant with respect to
the action by $\Gamma_1(3)$ on both sides.
\end{thm-int}

Here the $\Gamma_1(3)$-action on $\Stab^\dag(\DD_0)$ is induced by the subgroup $\Gamma_1(3) \subset
\Aut^\dag(\DD_0)$ identified in Theorem \ref{thmi:auto}.

On the level of central charges, the embedding is given in terms of a Picard-Fuchs differential
equation: for a fixed $E \in \DD_0$, the function $(\ZZ \circ I)(z)(E) \colon \tildeMG13 \to \C$ is a solution
of the Picard-Fuchs equation.
In particular, while classical enumerative mirror symmetry gives an interpretation of formal expansions of solutions of Picard-Fuchs equations at special
points of $M$ in terms of genus-zero Gromov-Witten invariants on $Y$, the space of stability conditions allows us to
interpret these solutions globally.

\subsection{Relation to existing work}
Various examples of stability conditions in local Calabi-Yau situations have been studied in the literature.
In particular, the local derived category of curves inside surfaces has been studied in
\cite{Thomas:stability,Bridgeland:ADE,Ishii-Uehara:An,Ishii-Ueda-Uehara,Okada:CY2,MMS:inducing,Brav-HThomas:ADE}, and results
similar to Theorem \ref{thmi:translates}, Theorem \ref{thmi:sc}, and Theorem \ref{thmi:auto} have been obtained.
Some examples of stability conditions on projective spaces were studied in \cite{Macri:Curves,Aaron-Daniele,Ohkawa}.
Other local Calabi-Yau threefold cases were studied in \cite{Toda:stab-crepant_res, Toda:CY-fibrations}, and, as already mentioned,
an open subset of $\Stab^\dag(\DD_0)$ has been described in \cite{Bridgeland:stab-CY}.

However, our approach follows the ideas in \cite{Bridgeland:K3} more closely than most of the above mentioned articles,
as we describe stability conditions in terms of stability of sheaves on $\P^2$, rather than in terms of exceptional
collections and quivers.
Applying this approach in our situation is possible due to the classical results of
Dr\'ezet and Le Potier \cite{Drezet-LePotier}; in particular, the fractal boundary of $\ZZ(U)$ discussed in 
Section \ref{secti:geomstab} is directly due to their results. At the same time, Sections \ref{sec:AlgebraicStability1} and \ref{sec:AlgebraicStability2} rely 
heavily on the work in \cite{Goro-Ruda:Exceptional} on exceptional collection and mutations.

Stability conditions around the orbifold point can be understood in terms of stability of quiver representations as
studied in \cite{Alastair-Ishii}; in particular our Theorem \ref{thmi:translates} could be understood as a derived
version of \cite[Theorem 1.2]{Alastair-Ishii} applied to our situation. 

There does not seem to be an equivalent of Theorem \ref{thmi:MS} in the literature for a Calabi-Yau 3-category; however, it is motivated by the
conjectural picture described in \cite[Section 7]{Bridgeland:spaces}.

There are many articles in the mathematical physics literature related to $\Pi$-stability and mirror symmetry for
$\C^3/\Z_3$ and the local $\P^2$ (as well as other local del Pezzo surfaces), and our presentation in Section
\ref{sec:MS} is very much guided by \cite{Aspinwall:Dbranes-CY} and \cite{ABK:topological_strings}.  In particular,
Theorem \ref{thmi:MS} is based on the computations of analytic continuations and monodromy for solutions of the
Picard-Fuchs equation of the mirror of the local $\P^2$ in \cite{AGM:measuring, Aspinwall:Dbranes-CY,
ABK:topological_strings}; in some sense, we are just lifting their results from the level of central charges to the
level of stability conditions.

In order for this to work, the ``central charges predicted by physicists'' had to survive a non-trivial test: they had
to satisfy the inequalities of Definition \ref{def:setG} (see Observation (\ref{obs:inequality}), page
\pageref{obs:inequality}).  The fact that they survived this test is somewhat reassuring for the case
of compact Calabi-Yau threefolds: identifying similar inequalities (which would be based on
inequalities for Chern classes of stable objects), and checking that the
central charges satisfy them, is the major obstacle towards constructing stability conditions on compact Calabi-Yau
threefolds.

\subsection{Open questions}
Bridgeland's conjecture \cite[Conj.\ 1.3]{Bridgeland:stab-CY} remains open; it
would identify $\Stab^\dag(\DD_0)$ with an open subset of the extended Frobenius manifold of the quantum cohomology of
$\P^2$.  Theorems \ref{thmi:translates} and \ref{thmi:sc} of this paper essentially complete the study of
$\Stab(\DD_0)$ as started in \emph{loc.\ cit.}; and Theorem \ref{thmi:MS} clarifies the discussion in \emph{loc.\
cit.} about the ``small quantum cohomology locus'', as this locus corresponds to the image of $\tildeMG13$. What is
missing from a proof of the whole conjecture, as pointed out in \emph{loc.\ cit.}, is still a better understanding of
the Frobenius manifold side.

It seems natural to conjecture that the full group $\Aut(\DD_0)$ of autoequivalences of $\DD_0$ preserves the connected
component $\Stab^\dag(\DD_0)$; in fact, this last one may be the only three-dimensional component of
$\Stab(\DD_0)$.  In this case, Theorem \ref{thmi:auto} would give a complete description of $\Aut(\DD_0)$.

Maybe the most intriguing question about $\Stab^\dag(\DD_0)$ related to our results is whether there is an intrinsic 
characterization of the image of the map $I$ of Theorem \ref{thmi:MS}, a question raised in other contexts in
\cite{Bridgeland:spaces}. To this end, note that the central charge on the image can also be given in terms of 
an analytic continuation of the genus
zero Gromov-Witten potential of $X$ (see \cite{ABK:topological_strings, Iritani:survey}; that this agrees with our
description using the mirror is classical enumerative mirror symmetry). But the genus-zero Gromov-Witten
potential is in turn determined by counting invariants of one-dimensional torsion sheaves (\cite{PT1, Toda:generating}),
i.e., counting invariants of stable objects close to the large-volume limit.

It would also be interesting to generalize some of the results of this paper to other ``local del Pezzo surfaces''.
In such a case, the starting point would be a generalization of the result in \cite{Drezet-LePotier} on the description of the Chern classes of stable sheaves.
Already for $\P^1\times\P^1$ the situation is more complicate: see \cite{Rudakov:Quadric,Rudakov:DelPezzoSurfaces} for results in this direction.

\subsection{Plan of the paper}
The paper is organized as follows.  In Section \ref{sec:GeoStability} we define geometric stability conditions and state
Theorem \ref{thm:geom-stability}, which classifies them.  Sections \ref{sec:Constraining} and \ref{sec:constructing} are
devoted to the proof of Theorem \ref{thm:geom-stability}.  In Section \ref{sec:boundary} we describe the boundary $\partial U$ of the
geometric chamber and prove Theorem \ref{thmi:translates}.  Algebraic stability conditions are
introduced in Section \ref{sec:AlgebraicStability1} in order to prove Theorem \ref{thmi:sc} (whose proof will take Section \ref{sec:AlgebraicStability2}).

In Section \ref{sec:autoequivalences} we study the group of autoequivalences and
prove Theorem \ref{thmi:auto}. Section \ref{sec:MS} discusses how the previous results
fit into expectations from mirror symmetry, and includes the proof of Theorem \ref{thmi:MS}.
Finally, three appendices complete the paper.
In Appendix \ref{app:DP}, we review the results of Dr\'ezet and Le Potier
as we need them in the proof of Theorem \ref{thm:geom-stability}.
Appendix \ref{app:BridgelandFramework} is a brief
introduction to stability conditions and contains an improved criterion for the existence of Harder-Narasimhan
filtrations.
In Appendix \ref{app:ineq} we give a sketch of the proof that the central charges we define in Section \ref{sec:MS} satisfy the inequalities in Definition \ref{def:setG}.

\subsection{Notation}
We work over the complex numbers $\C$.
We let $X$ denote the total space of $\OO_{\P^2}(-3)$, and
$i\colon\P^2\into X$ the inclusion of the zero-section.
We let
$\Coh_0:=\Coh_{\P^2}X\subset\Coh X$ be the subcategory of coherent sheaves on $X$ supported
(set-theoretically) on the zero-section. We write 
$\DD_0 = \Db_0(X)$ for the subcategory of $\Db(\Coh X)$ of complexes
with bounded cohomology, such that all of its cohomology sheaves
are in $\Coh_0$. (Note that $\DD_0 \cong \Db(\Coh_0)$ as observed in \cite{Ishii-Uehara:An}, Notation and Convention.)
The space of stability conditions on $\DD_0$ will be denoted by $\Stab(\DD_0)$, and its two subsets of geometric and
algebraic stability conditions by $U$ and $\Stab_a$, respectively (see Definitions \ref{def:GeomStability} and
\ref{def:algebraicstability}).

An object $S$ in $\DD_0$ is called \emph{spherical} if $\Ext^p(S,S)\cong\C$ for $p=0,3$ and is zero otherwise.
For a spherical object $S$ we denote by ${\ST}_S$ the spherical twist associated to $S$, defined by the exact triangle
\[
\Hom^*(S,M)\otimes S\stackrel{ev}{\longrightarrow} M\longrightarrow {\ST}_S(M),
\]
for $M\in\DD_0$ (see \cite{Seidel-Thomas:braid}).

By abuse of notation, we will write
$\OO_{\P^2}(n) \in \DD_0$ for the spherical objects $i_*\OO_{\P^2}(n)$.
For $x \in \P^2$ we denote by $k(x)$ the skyscraper sheaf in $X$ of length one
concentrated at $x$.

The Grothendieck group of $\DD_0$ is denoted by $K(\DD_0)$. It is isomorphic to $\Z^{\oplus 3}$.
For any $E \in\DD_0$, we write $r(E), d(E), c(E)$ for the components of
the Chern character of its push-forward to $\P^2$; more precisely, if $\pi\colon X\to\P^2$ is the
projection and $\ch \colon K(\P^2) \to A_*(\P^2) \otimes \Q$ the Chern character with values in the Chow ring, then we write
\[
\ch(\pi_*(E)) = r(E) \cdot [\P^2] + d(E) \cdot [\mathrm{line}]+ c(E) \cdot \pt.
\]

For a complex number $z\in\C$, we write $\Re z$ (resp.\ $\Im z$) for its real (resp.\ imaginary) part.

\subsection{Acknowledgements}
It is a pleasure to thank Hiroshi Iritani, from whom we first learned about the possibility of a $\Gamma_1(3)$-action on
$\DD_0$,  and Aaron Bertram, Tom Bridgeland, Sukhendu Mehrotra, Paolo Stellari, and Richard Thomas for very interesting and
useful discussions. We would also like to thank Eric Miles and Ryo Ohkawa for valuable comments on a preliminary
version of the paper, and the referees for their detailed and very helpful
comments which greatly improved the manuscript and corrected many inaccuracies.

We are very grateful to the Mathematical Sciences Research Institute in Berkeley, and its
Jumbo program in algebraic geometry in the Spring 2009, as several of the discussions mentioned
above took place during the first author's stay at the MSRI. We would also like to thank
the University of Bonn for
the warm hospitality during the writing of parts of this paper. The first author was partially supported by the NSF
grant DMS-0801356/DMS-1001056 and the second author by the NSF grant DMS-1001482, the Hausdorff Center for Mathematics, Bonn, and SFB/TR 45.

\section{Geometric stability conditions}\label{sec:GeoStability}

We assume familiarity with the
notion of stability conditions on a derived category; see Appendix \ref{app:BridgelandFramework} for a short summary, and
\cite{Bridgeland:Stab}, \cite[Section 3.4]{Kontsevich-Soibelman:stability} for a complete reference.

We begin by constructing and classifying ``geometric'' stability conditions on $\DD_0$.
Loosely speaking, geometric stability conditions are those that are most closely connected to the geometry of 
sheaves on $X$; in the definition below, we require that the simple objects of $\Coh_0$ remain stable, but it will
also turn out that the semistable objects are at most two-term complexes of sheaves in $\Coh_0$.

\begin{Def}\label{def:GeomStability}
A stability condition $\sigma$ on $\DD_0$ is called \emph{geometric} if the
following two conditions are satisfied:
\begin{enumerate}
\item All skyscraper sheaves $k(x)$ of closed points $x\in\P^2$ are
$\sigma$-stable of the same phase.
\item \label{enum:full}
The connected component of $\Stab(\DD_0)$ containing $\sigma$ is full, that is it has maximal dimension equal to $3$.
\end{enumerate}
\end{Def}

We write $U$ for the set of geometric stability conditions, and refer to it as the ``geometric chamber''; in
fact, we will see that it is precisely one of the chambers with respect to the chamber decomposition given by the
wall-crossing phenomenon for semistable objects of class $[k(x)]$.

Part (\ref{enum:full}) of the definition is a technical condition to ensure that the wall-crossing for semistable
objects behaves nicely (see \cite[Section 9]{Bridgeland:K3}); it is equivalent to the ``support
property'' introduced by Kontsevich-Soibelman (see Proposition \ref{prop:SupportProperty}).

We recall that a Bridgeland stability condition can be constructed by giving the heart of a bounded t-structure
$\AA \subset \DD_0$, and a compatible central charge $Z \colon K(\AA) = K(\DD_0) \to \C$ that sends objects
in $\AA$ to the semi-closed upper half plane (see Remark \ref{rmk:tstruct}). The t-structures appearing in geometric
stability conditions are given by the now familiar notion of \emph{tilting} (see \cite{Happel-al:tilting}):

For purely 2-dimensional sheaves $\FF \in \Coh_0$, the slope
function $\mu(\FF) = \frac{d(\FF)}{r(\FF)}$ gives a notion of slope-stability (as in Definition
\ref{def:stablesheaf}).
By the same arguments as in the case of a projective variety, Harder-Narasimhan filtrations exist.
Thus we can follow \cite[Lemma 5.1]{Bridgeland:K3} to make the following definition:

\begin{Def}
For any $B \in \R$, let
$\left(\Coh_0^{>B}, \Coh_0^{\le B}\right)$ be the torsion pair in $\Coh_0$ determined by:
\begin{itemize}
\item $\Coh_0^{\le B}$ is generated (by extensions) by semistable sheaves of slope
$\mu \le B$, and
\item $\Coh_0^{>B}$ is generated by semistable sheaves of slope
$\mu > B$ and zero- or one-dimensional torsion sheaves.
\end{itemize}
Let $\Coh_0^{\sharp(B)} \subset \DD_0$ be the tilt of $\Coh_0$ at the
torsion pair $\left(\Coh_0^{>B}, \Coh_0^{\le B}\right)$, that is
\[
\Coh_0^{\sharp(B)}=\left\{E\in\DD_0\,\colon\,\begin{array}{l}\bullet\ \HH^i(E)=0,\text{ for all }i\neq0,-1\\ \bullet\ \HH^0(E)\in\Coh_0^{>B}\\ \bullet\ \HH^{-1}(E)\in\Coh_0^{\le B}\end{array}\right\}.
\]
\end{Def}

The structure of central charges compatible with $\Coh_0^{\sharp(B)}$ (and thus, as we will see,
the structure of the whole geometric chamber) depends tightly on the set
of Chern classes for which there exist stable torsion-free sheaves.
In the case of $\P^2$, Dr\'ezet and Le Potier have given a complete description of this
set (see Appendix \ref{app:DP} for more details).
It is most naturally described in terms of the discriminant $\Delta(\FF)$, which is defined as
\[
\Delta(\FF) = \frac{d(\FF)^2}{2r(\FF)^2} - \frac{c(\FF)}{r(\FF)}.
\]
Recall that a vector bundle $\EE$ on $\P^2$ is called \emph{exceptional} if
$\Hom(\EE, \EE) = \C$, and $\Ext^p(\EE, \EE) = 0$ for $p > 0$. 
For an exceptional vector bundle $\EE_\alpha$ of rank $r_{\alpha}$ and slope $\alpha$, it follows from
Riemann-Roch that the discriminant is given by $\Delta_\alpha = \frac 12 - \frac 1{2r_\alpha^2}$. 
Also note that an exceptional vector bundle can equivalently be characterized by being slope-stable with
discriminant smaller than $1/2$: see \cite[Theorem 4.1]{Goro-Ruda:Exceptional} and \cite[Lemme (4.2)]{Drezet-LePotier}. 
As the slopes of exceptional vector bundles can be constructed explicitly (see Theorem \ref{thm:DP-A}), it remains to describe the slopes and discriminants
of non-exceptional stable torsion-free sheaves.

As explained in \cite[Section 16]{LePotier}, one can slightly reformulate the results of
\cite{Drezet-LePotier} and construct
a function $\delta_\infty^{DP} \colon \R \to [\frac 12, 1]$. It is periodic of period $1$ and
Lipschitz-continuous with Lipschitz constant $\frac 32$.
We refer to Appendix \ref{app:DP} for the precise definition of $\delta_\infty^{DP}$.
Its construction is motivated by the following observation: If $\FF_\beta$ is a slope-stable
sheaf with $\beta < \alpha$, then $\Hom(\EE_\alpha, \FF_\beta) = 0$; if we additionally assume
$\beta > \alpha - 3$, then, by Serre duality, also $\Ext^2(\EE_\alpha, \FF_\beta) = 0$ and hence
$\chi(\EE_\alpha, \FF_\beta) \le 0$. Using Riemann-Roch this yields an inequality of the form
$\Delta(\FF_\beta) \ge p_\alpha(\beta)$ for $\alpha - 3 < \beta < \alpha$, with $p_\alpha(x)$ being a quadratic
polynomial.  The function $\delta_\infty^{DP}$ is the supremum of all the quadratic polynomials $p_\alpha$
restricted to the ranges where the inequality is valid. 

The main result of \cite{Drezet-LePotier} is that $\Delta \ge \delta_\infty^{DP}(\mu)$ is not only a necessary,
but also a sufficient condition for the existence of a stable torsion-free sheaf of slope $\mu$ and discriminant
$\Delta$. For later use, we paraphrase their result as follows:

Define $S_\infty \subset \R^2$ to be the closed subset lying above the
the graph of $\delta_\infty^{DP}$, i.e.,
\[
S_\infty = \stv{(\mu, \Delta) \in \R^2}{\Delta \ge \delta_\infty^{DP}(\mu)}.
\]
\begin{Thm}[\cite{Drezet-LePotier}] \label{thm:DP}
Let $S \subset \Q^2$ be the set of pairs $(\mu(\FF), \Delta(\FF))$ 
where $\FF$ is any slope-stable torsion-free sheaf $\FF$ on $\P^2$.
Similarly, let $S_E \subset \Q^2$ be the corresponding set for
slopes and discriminants of 
exceptional vector bundles.  Then $S$ is the disjoint union
\[
S = S_E \dotcup \left( S_\infty \cap \Q^2 \right)
\]
The set $S_E$ has no accumulation points in $\R^2 \setminus S_\infty$.
\end{Thm}

\begin{figure}[htb]
\begin{center}
\includegraphics{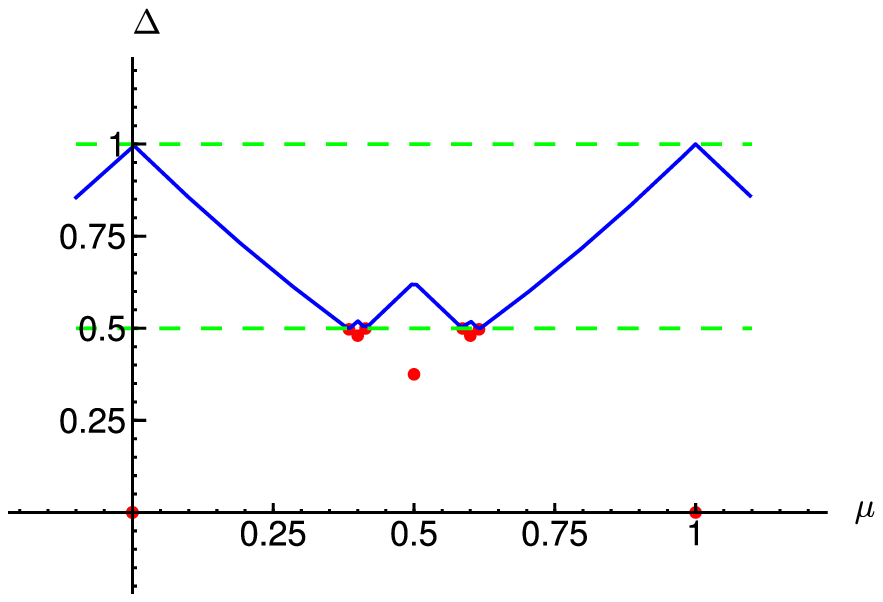}
\caption{$\delta_\infty^{DP}$ and exceptional objects}
\label{fig:DP-plot}
\end{center}
\end{figure}

We explain this reformulation of Dr\'ezet and Le Potier's result
in Appendix \ref{app:DP}, along with their explicit description of the set $S_E$; see also Figure \ref{fig:DP-plot}.

\begin{Def}\label{def:setG}
We define the set $G \subset \C^2$ as the set of pairs $a, b \in \C$ satisfying the following three
inequalities (where we set $B := -\frac{\Im b}{\Im a}$ if $a$ satisfies the first inequality):
\begin{align}
\Im a &> 0, 			\label{ineq:Ima}\\
\Re b &> -B \cdot \Re a - \delta_\infty^{DP}(B) + \frac 12 B^2,     \label{ineq:genslope}
\intertext{and, in case there exists an exceptional vector bundle of slope
$B$ and discriminant $\Delta_B$,}
\Re b &> -B \cdot \Re a - \Delta_B + \frac 12 B^2.    \label{ineq:exceptslope}
\end{align}
\end{Def}

\begin{Thm}\label{thm:geom-stability}
For $a, b \in \C$, denote by $Z_{a, b} \colon K(\DD_0) \to \C$ the central charge given by
\begin{equation}\label{eq:Zab}
Z_{a,b}(E)=-c(E)+ad(e)+br(E)
\end{equation}
for $E \in \DD_0$.

Then there exists a geometric stability condition $\sigma_{a,b} = (Z_{a,b}, \PP_{a,b})$ with
$Z_{a, b}$ as above if and only if $(a, b) \in G \subset \C^2$.  Its heart is, up to shifts, given by
$\PP_{a,b}((0,1]) = \Coh_0^{\sharp(B)}$.  The shifts $k(x)[n]$ of skyscraper sheaves are the only
stable objects of class $\pm [k(x)]$.

Any geometric stability condition is equivalent to a stability condition $\sigma_{a,b}$ up to
the action of a unique element in $\C$.
\end{Thm}

The action of $z \in \C$ is given in Remark \ref{rmk:GroupAction}: it is the lift to the space of stability conditions of the
multiplication by $\exp(z)$ on $\Hom(K(\DD_0), \C)$.

The theorem can be rephrased as stating that
$U/\C \cong G$, with a section given by $(a, b) \mapsto \sigma_{a,b}$.
Later, in Remark \ref{Rem:sc-slice}, we will see that this slice of the $\C$-action can be extended to a whole connected
component of $\Stab(\DD_0)$.  The theorem will be proved in the following two sections.

The best visualization of the set of allowed central charges is given by the
following observation: as long as $\Im a > 0$, the central charge can be
thought of as a surjective map $K(\DD_0)_\R \cong \R^3 \to \C \cong
\R^2$.  Up to the action of $\GL_2^+(\R)$ on $\R^2$ (which does not affect the
set of stable objects), this map is determined by its kernel, and by the
orientation induced on $K(\DD_0)_\R / \Ker Z$. As long as $\Im a > 0$,
the orientation does not change.
The kernel intersects the affine hyperplane $r = 1$ of
$K(\DD_0)_\R$ in a single point. 
The inequalities are equivalent to requiring that this point lies below
the graph of $\Delta = \delta_\infty^{DP}(\mu)$, and not on any of the
rays going up vertically from a point in $S_E$ (see Figure \ref{fig:DP-plot}).

However, for several reasons it is helpful to classify geometric stability conditions up to the
action of $\C$, rather than the action of $\GL_2^+(\R)$. The subgroup $\C$ acts on $\Stab(\DD_0)$ with closed orbits, it has no stabilizers, and it has a well-behaved quotient
$\Stab(\DD_0)/\C$. The stability conditions $\sigma_{a, b}$ constitute a slice of this action on
$U \subset \Stab(\DD_0)$, and the boundary of $U$ can be identified, up to
the $\C$-action, with the boundary of the set $\stv{\sigma_{a, b}}{(a,b) \in G}$. None of these
statements would hold for the $\GL_2^+(\R)$-action, and the picture of the preceding paragraph only
gives a partial picture of the boundary of $U$: one can see
that every ray starting at a point of $S_E$ going up vertically may give two walls in the boundary
$\partial U$ of the set $U$, but we cannot see that many of these walls intersect at
points where the central charge lies in the real line.

\section{Constraining geometric stability conditions}\label{sec:Constraining}

In this section we will show that geometric stability conditions can only be
of the form given in Theorem \ref{thm:geom-stability}. The general idea is the same as in \cite{Bridgeland:K3}: 
if we assume that the skyscraper sheaves $k(x)$ are stable of phase 1, then $\Hom$-vanishing helps to constrain the form
of objects in $\PP((0,1])$, and we can identify $\PP((0,1])$ with an explicit tilt of the standard t-structure. By the
existence of a well-behaved chamber decomposition for the wall-crossing for stable objects of class $[k(x)]$, the set of
geometric stability conditions is open, and we need to prove inequalities for the central charge only when it
is defined over $\Q$.

The proof will be broken into several lemmata and propositions. The following
observation shows that the bound of Theorem \ref{thm:DP} translates into
bounds for stable objects in $\Coh_0$:
\begin{Lem}\label{lem:mustabilitypf}
A sheaf $\FF \in \Coh_0$ is a pure slope-stable sheaf
if and only if it is the push-forward $\FF = i_* \FF_0$ of some
slope-stable pure sheaf $\FF_0 \in \Coh \P^2$.
\end{Lem} 
\begin{Prf}
Since $i_* \colon \Coh \P^2 \into \Coh_0$ is
a full subcategory, closed under subobjects and quotients, and
since $i_*$ preserves the ordering by slopes, it follows $\FF_0 \in \Coh \P^2$
is stable if and only if $i_* \FF_0$ is stable.

Now assume that $\FF$ is stable. Then $\End \FF = \C \cdot \Id$.
Let $Z$ be the scheme-theoretic support of $\FF$.
By definition, its global sections act faithfully on $\FF$, so $H^0(\OO_Z) \cong \C$.
Hence $Z$ must be contained scheme-theoretically in the fiber of the origin
under the contraction $X \onto \Spec H^0(\OO_X) \cong \C^3/\Z_3$, as otherwise the image
of $H^0(\OO_X) \to H^0(\OO_Z)$ would be non-trivial.
But the scheme-theoretic fiber of the origin is exactly $\P^2$, and so $\FF$ is the push-forward 
$i_* \FF_0$ of some sheaf $\FF \in \Coh \P^2$ on $\P^2$.
\end{Prf}

Now assume we are given a geometric stability condition.  After a rescaling by
$\C$, we may assume that all skyscraper sheaves $k(x)$ of closed points $x\in\P^2$
are stable with phase $1$ and $Z(k(x)) = -1$.

\begin{Lem}[{\cite[Lemma 10.1]{Bridgeland:K3}}]\label{lem:geom-stability} 
Let $(Z,\PP)$ be a stability condition such that the skyscraper sheaves
$k(x)$ are stable of phase 1, with $Z(k(x)) = -1$.
\begin{enumerate}
\item\label{coh-constraint}
For any object $E \in \PP((0,1])$, its cohomology sheaves $\HH^i(E)$
vanish unless $i = 0, -1$.
\item 					\label{H-1_pure}
Further, for any such $E \in \PP((0,1])$ the cohomology 
sheaf $\HH^{-1}(E)$ is pure of dimension 2.
\item 					\label{stableinP1}
If $E \in \PP(1)$ is stable and $E \neq k(x)$ for all $x \in \P^2$, then
there is a vector bundle $\FF$ on $\P^2$ such that $E \cong i_* \FF[1]$.
\end{enumerate}
\end{Lem}

\begin{Prf}
If $E \in \PP((0,1))$, then $\Hom(E, k(x)[i]) = 0$ for $i < 0$, and
$\Hom(k(x)[i], E) = \Hom(E, k(x)[3+i]) = 0$ for $i \ge 0$ and $x \in \P^2$.
Since $E$ is supported on $\P^2$, all homomorphisms with shifts of skyscraper sheaves outside the zero-section are zero.
We can therefore apply \cite[Prop.\ 5.4]{Bridgeland-Maciocia:K3Fibrations} and deduce that $E$ is quasi-isomorphic to a 3-term complex of
locally free sheaves $E^{-2}\stackrel{d^{-2}}{\longrightarrow} E^{-1} \to E^0$. Hence
$\HH^{-2}(E)$ torsion-free on $X$; since $\HH^{-2}(E) \in \Coh_0$,
it must vanish. This implies the first claim for such $E$.

This also shows that $\HH^{-1}(E)$ is the cokernel of an injective map
from a locally free sheaf to a torsion-free sheaf, which implies the second claim.

If $E$ is stable of phase 1, then additionally $\Hom(E, k(x)) = 0$ (as
they are both stable objects of the same phase).
Hence $E$ is isomorphic to a two-term complex of vector bundles
$E^{-2}\stackrel{d^{-2}}{\longrightarrow}E^{-1}$. By the same argument as in the previous case,
the map $d^{-2}$ must be injective, so that $E$ is isomorphic to the shift
of a sheaf: $E \cong \FF'[1]$. Since $E$ is stable, $\FF'$ can only have
scalar endomorphisms, and thus $\FF'$ is the push-forward of a sheaf
$\FF$ on $\P^2$. Using 
$0 = \Hom(i_* k(x),i_* \FF[1]) \cong \Hom(k(x) \oplus k(x)[1],\FF[1])$,
it follows that $\FF$ is a vector bundle.

Since all the assertions of the lemma are properties that are closed under
extensions, this finishes its proof.
\end{Prf}

Recall that, for a full stability condition $\sigma=(Z,\PP)$,
\begin{equation} 				\label{eq:def-metric}
\| W\|_{\sigma}:=\sup\left\{\frac{|W(E)|}{|Z(E)|}\colon E\text{ is }\sigma\text{-stable}\right\}
\end{equation}
defines a metric on $\Hom(K(\DD), \C)$ (see Appendix \ref{app:BridgelandFramework} for more
details).

The next result is based on \cite[Section 9]{Bridgeland:K3}:
\begin{Prop}\label{prop:chambers}
Let $\DD$ be a triangulated category such that $K(\DD)$ is a finite-dimensional
lattice, and let $\Stab^* \subset \Stab(\DD)$ be a full connected component of its space of
stability conditions.
Fix a primitive class $\alpha \in K(\DD)$, and an arbitrary set $S \subset \DD$ of objects of class
$\alpha$.  Then there exists a collection
of walls $W^S_\beta$, $\beta \in K(\DD)$, with the following properties:
\begin{enumerate}
\item Every wall $W^S_\beta$ is a closed submanifold with boundary of real
codimension one.
\item The collection $W^S_\beta$ is locally finite (i.e., every compact subset $K \subset \Stab^*$
intersects only a finite number of walls).
\item \label{enum:semistable-on-wall}
For every stability conditions $(Z, \PP) \in W^S_\beta$, there exists a phase $\phi$
and an inclusion $F_\beta \into E_\alpha$ in $\PP(\phi)$ with $[F_\beta] = \beta$ and some $E_\alpha \in S$.
\item \label{enum:chambers}
If $C \subset \Stab^*$ is a connected component of the complement of
$\bigcup_{\beta \in K(\DD)} W^S_\beta$, and $\sigma_1, \sigma_2 \in C$, then an object $E_\alpha \in S$
is $\sigma_1$-stable if and only if it is $\sigma_2$-stable. 
\end{enumerate}
\end{Prop}
\begin{Prf}
For a class $\beta \in K(\DD)$ let $V^S_\beta$ be the set of stability conditions $(Z, \PP)$
for which there exists an inclusion as in part (\ref{enum:semistable-on-wall}). Since $\alpha$ is
primitive, each $V^S_\beta$ is contained in the codimension-one subset with
$\Im \frac{Z(\beta)}{Z(\alpha)} = 0$.

We first want to show that there are only finitely many $\beta$ for which $V^S_\beta$ intersects an open
ball $B_{\frac 18}(\sigma)$ of radius $\frac 18$ around $\sigma = (Z, \PP)$: Given $\sigma$ and $S$,
let $I_\sigma(S) \subset K(\DD)$ be the set of all classes $\beta$ for which there exists $\phi \in
\R$ with $Z(\alpha) \in \R_{>0}\cdot e^{i\pi \phi}$ and a strict inclusion $F_\beta \into E$ in the
quasi-abelian category $\PP((\phi - \frac 14, \phi + \frac 14))$ with $[F_\beta]=\beta$ and $E \in
S$.
Since the metric $\|\cdot\|_\sigma$ is finite, and since $K(\DD)$ is a discrete
subgroup of $K(\DD)\otimes \R \cong \R^n$, there exist only finitely many classes $\gamma \in
K(\DD)$ that have a $\sigma$-semistable object $F_\gamma$ of class $[F_\gamma] = \gamma$ satisfying
$\abs{Z(F_\gamma)} < \abs{Z(\alpha)}$. It follows that the set $I_\sigma(S)$ is also finite (as each
HN filtration factor of $F_\beta$ is an object $F_\gamma$ as considered in the previous sentence).
But if $V^S_\beta$ intersects $B_{\frac 18}(\sigma)$,
then it follows from \cite[Lemma 7.5]{Bridgeland:Stab} that $\beta \in I_\sigma(S)$.

An object $E$ of class $\alpha$ is $(Z', \PP')$-semistable for $(Z', \PP') \in B_{\frac 18}(\sigma)$ if and only
if $\Im \frac{Z'(\beta)}{Z'(\alpha)} \le 0$ for 
every $\beta \in I_\sigma(\{E\})$---and it is stable if and only if the inequalities are strict.
Repeating this argument for every possible subobject $F_\beta$, it follows that inside the codimension one
subset $\Im \frac{Z'(\beta)}{Z'(\alpha)} = 0$, the set $V^S_\beta$ is a finite union of subsets, each
of which is  cut out by a finite number of inequalities of the form $\Im \frac{Z'(\beta')}{Z'(\alpha)}
\le 0$ for some $\beta' \in I_\sigma(S)$.  We let $W^S_\beta$ be the union of all
codimension-one components of $V^S_\beta$.

It remains to prove claim (\ref{enum:chambers}). It is sufficient to consider the case $\sigma_1, \sigma_2 \in
B_{\frac 18}(\sigma) \cap C$. Assume that there is an object $E \in S$ that is $\sigma_1$-stable but
not $\sigma_2$-stable.  Then on every path $\gamma \colon [0,1] \to B_{\frac 18}(\sigma) \cap C$
connecting $\sigma_1$ with $\sigma_2$, there is a point $\gamma(t)$ on which $E$ is strictly
semistable, i.e., $\gamma(t) \in V^S_\beta \cap C$ for some $\beta \in I_\sigma(S)$ and $t \in (0, 1]$.  But
by the definition of the walls $W^S_\beta$, the set $V^S_\beta \cap C$ has codimension at least two,
and hence we may choose $\gamma$ such that for $t \in (0, 1)$, it avoids all of the finitely many
non-empty subsets $V^S_\beta \cap C \subset C$ for $\beta \in I_\sigma(S)$, in other words we have
that $E$ is $\gamma(t)$-stable for $t \in (0, 1)$, and 
$\sigma_2 \in V^S_\beta \cap C$ for some $\beta \in I_\sigma(S)$. In particular, $\sigma_2$ is contained
in the set $\Im \frac{Z(\beta)}{Z(\alpha)} = 0$, and $E$ will not be stable in the subset
of $B_{\frac 18}(\sigma) \cap C$ with $\Im \frac{Z(\beta)}{Z(\alpha)} \le 0$. On the other hand,
the set $C \setminus \bigcup V^S_\beta$ is path-connected, and by the previous argument
$E$ is stable on all of it.  This is a contradiction.

(In other words, we proved in the last step that higher-codimension components of $V^S_\beta$ always
come from objects $E_\alpha$ that are semistable on this component, and unstable at any nearby point.)
\end{Prf}

As can be seen from the proof, the proposition holds for any family of stability conditions
satisfying the \emph{support property} of \cite[Section 1.2]{Kontsevich-Soibelman:stability}.
In fact, in our situation fullness and the support property are equivalent, see Proposition
\ref{prop:SupportProperty}.

\begin{Cor}\label{cor:openness}
The set $U$ of geometric stability conditions is open in the space of stability
conditions $\Stab(\DD_0)$. Its boundary $\partial U = \overline{U} \setminus U$
is given by a locally finite union of walls, and each wall is a real submanifold with boundary
in $\Stab(\DD_0)$ of codimension one.
\end{Cor}

We proceed to show that any geometric stability condition is necessarily of the form given in Theorem \ref{thm:geom-stability}.

By Lemma \ref{lem:geom-stability}, we have
$\PP((0,1]) \subset \langle \Coh_0, \Coh_0[1] \rangle$. 
This implies that $\PP((0,1])$ is obtained from $\Coh_0$ by tilting at the torsion pair
\begin{align*}
\TT &= \Coh_0 \cap \PP((0,1]) \\
\FF &= \Coh_0 \cap \PP((-1,0])
\end{align*}
(see e.g.\ \cite[Lemma 1.1.2]{Polishchuk:families-of-t-structures}).

Since we assume $Z(k(x)) = -1$, the central charge can be written in the
form of equation \eqref{eq:Zab}; in particular
$\Im Z(E) = d(E) \cdot \Im a + r(E) \Im b$.
By mimicking the proof of \cite[Prop.\ 10.3]{Bridgeland:K3}, it follows that $\Im a>0$ and, after setting $B=-\frac{\Im b}{\Im a}$, $\PP((0,1]) = \Coh_0^{\sharp(B)}$.

It remains to prove the inequality on $\Re b$.  We first assume that $B \in
\Q$. For any semistable torsion-free sheaf $\FF$ on $\P^2$ of slope $B$ we
have $\Im Z(i_*\FF) = 0$ and $i_*\FF[1] \in \Coh_0^{\sharp(B)}$, hence we must
have $\Re Z(i_*\FF) > 0$. It follows that:
\begin{align*}
0 & < \frac{\Re Z(i_*\FF)}r = \Re b + \Re a B - \frac cr \\
&= \Re b + \Re a B + \Delta(\FF) - \frac 12 B^2
\end{align*}
Applying Theorem \ref{thm:DP}, we obtain the inequalities
(\ref{ineq:genslope}) and (\ref{ineq:exceptslope}).

Finally, we need to treat the case $B \not \in \Q$. By Corollary
\ref{cor:openness}, there exists an open neighborhood
$V \subset \C^2$ of $(a,b)$, such
that any $(a', b') \in V$ with
$\frac{\Im b'}{\Im a'} \in \Q$ satisfy inequality (\ref{ineq:genslope}).
Hence it holds for $(a, b)$, too.

\section{Constructing geometric stability conditions}\label{sec:constructing}

We now come to the proof of existence of geometric stability conditions. The main problem is to prove the existence of
Harder-Narasimhan filtrations for the stability function $Z_{a, b}$ on $\Coh_0^{\sharp(B)}$. We
prove this directly in the case where the image of $\Im Z$ is discrete, and then use Bridgeland's deformation result
to extend it to the more general case. In order to make the extension effective, we have to bound the metric
$\|\cdot\|_\sigma$ on $\Hom(K(\DD_0), \C) \cong \C^3$ defined by equation \eqref{eq:def-metric} from above.
To do so, we in turn have to control $\abs{Z(E)}$ for stable objects $E$ from below.

Our arguments in this section build on \cite{Bridgeland:K3, Aaron-Daniele}. 

\begin{figure}[htb] \begin{center}
	\includegraphics{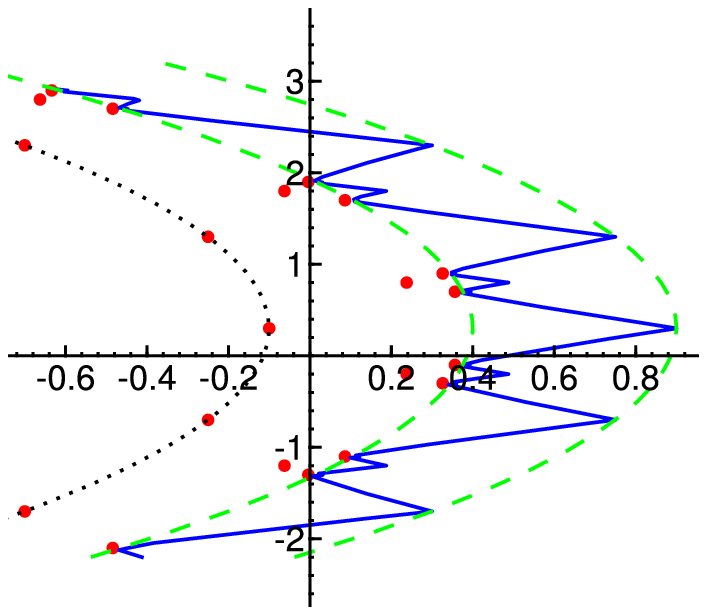}
\caption{The path $\gamma_{a,b}(t)$ (blue) and the central charges $\frac{Z(E)}{r(E)}$ of
exceptional vector bundles $E$ (red)}
	\label{fig:pathgamma}
\end{center} \end{figure}

Given $(a, b) \in G \subset \C^2$, let $\gamma_{a,b}\colon\R\to\C$ be the infinite path $\gamma_{a,b}(t)=x(t)+iy(t)$ defined by
\begin{equation*}\label{eq:pathgamma}
\begin{split}
x(t) &= \Re b + \frac 12 (\Re a)^2 + \delta_\infty^{DP}(t)
	-\frac 12 \left(t - \Re a\right)^2 \\
y(t) &= \Im a\cdot t + \Im b
\end{split}
\end{equation*}
and let $S_{a, b}\subset\C$ be the closed subset cut out by $\gamma_{a,b}$
that lies on or to the right of $\gamma_{a,b}$, i.e.,
\[
S_{a, b}:=\stv{x + iy}{\exists t\text{ with }y=y(t),\, x\ge x(t)}.
\]
Since $\delta_\infty^{DP}(t) \in [\frac 12, 1]$, the path $\gamma_{a,b}$
is contained between the graphs of two parabolas with horizontal distance
$\frac 12$, see Figure \ref{fig:pathgamma}.

\begin{Lem}\label{lem:ineq}
Let $\FF\in\Coh \P^2$ be a torsion-free slope-stable sheaf on $\P^2$
of rank $r$ that is not an exceptional vector bundle.
Then $\frac{Z_{a,b}(i_*\FF)}r \in S_{a,b}$.
\end{Lem}

\begin{Prf}
We write $\frac {\ch(\FF)}r = [\P^2] + \mu [l] +  \frac cr \pt$
and $\frac{Z_{a,b}(i_*\FF)}r = x + iy$. Using
Theorem \ref{thm:DP}, we obtain:
\begin{align*}
y &= \Im a \mu + \Im b \\
x &= - \frac cr + \Re a \mu + \Re b  = \Delta(\FF) - \frac 12 \mu^2 + \Re a \mu + \Re b \\
& \ge \delta_\infty^{DP}(\mu) - \frac 12 \mu^2 + \Re a \mu + \Re b
\end{align*}
Setting $t = \mu$ yields the claim.
\end{Prf}

\begin{Lem}
For any $(a, b) \in G$, the central charge $Z_{a,b}$ is a stability function for
$\Coh_0^{\sharp(B)}$: if $0\neq E\in\Coh_0^{\sharp(B)}$, then $Z_{a,b}(E)\in\H$.
\end{Lem}
\begin{Prf}
It is sufficient to prove the claim for torsion sheaves of dimension $\le 1$, and for (shifts of) slope-stable sheaves.
The claim is only non-trivial for objects $\FF[1]$, where $\FF$ is
a purely $2$-dimensional slope-stable of slope $\mu(\FF) = B$. In such a case, $Z_{a,b}(\FF)$
is lying on the real line. By Lemma \ref{lem:mustabilitypf}, $\FF$
is the push-forward of a slope-stable sheaf on $\P^2$. If $\FF$
is not an exceptional vector bundle, the previous lemma shows
$Z_{a,b}(i_*\FF) \in \R_{>0}$. If it is exceptional, the same computation
yields this statement from inequality \eqref{ineq:exceptslope}.
\end{Prf}

In other words, inequality \eqref{ineq:genslope} is equivalent to
$0 \notin S_{a,b}$, i.e., the path is passing through the real line with
positive real part; and inequality \eqref{ineq:exceptslope} is equivalent
to $Z(\EE) \not\in \R_{\le 0}$ for any exceptional vector bundle $\EE$ on
$\P^2$. Together, they guarantee that the central charge of a slope-stable
sheaf never lies on the negative real line.

Due to the above Lemma, as in Remark \ref{rmk:tstruct} we can define a phase function $\phi(E):=(1/\pi)\arg(Z(E))\in(0,1]$, for all $0\neq E\in\Coh_0^{\sharp(B)}$.

\begin{Lem}[{\cite[Proposition 7.1]{Bridgeland:K3}}] 		\label{lem:HNFiltrationsDiscrete} 
Let $a,b\in\C$ be such that $B\in\Q$ and $Z_{a,b}(\Coh_0^{\sharp(B)}\setminus\{0\})\subset\H$.
Then Harder-Narasimhan filtrations exist for $(Z_{a,b},\Coh_0^{\sharp(B)})$.
\end{Lem}
\begin{Prf}
We use Proposition \ref{prop:HNFiltrationsDiscrete}.
Since $B\in\Q$, the image of $\Im(Z_{a,b})$ is discrete in $\R$.

Let $F\in\Coh_0^{\sharp(B)}$ and let
\begin{equation}\label{eq:Bonn17}
0=A_0\subset A_1\subset\ldots\subset A_j\subset A_{j+1}\subset\ldots\subset F,
\end{equation}
be a sequence of subobjects, with $A_j\in\PP'_{a,b}(1)$. As in Proposition \ref{prop:HNFiltrationsDiscrete}, $\PP'_{a,b}(1)$ is the full subcategory of $\Coh_0^{\sharp(B)}$ whose objects have phase $1$ with respect to $Z_{a,b}$.
We need to show that \eqref{eq:Bonn17} stabilizes.

To this end, we first observe that \eqref{eq:Bonn17} induces a sequence of inclusions
\begin{equation}\label{eq:Bonn18}
0=\HH^{-1}(A_0)\into\ldots\HH^{-1}(A_j)\into\ldots\into\HH^{-1}(F).
\end{equation}
Since $\Coh_0$ is Noetherian, \eqref{eq:Bonn18} must terminate.
We can therefore assume that $\HH^{-1}(A_j)\cong\FF_{-1}$ for all $j$ and some $\FF_{-1}\in\Coh_0^{\le B}$.
Let $\FF_0$ denote the cokernel of the inclusion $\FF_{-1}\into\HH^{-1}(F)$; 
then, by the long exact cohomology sequence, we have $\FF_0\in\Coh_0^{\le B}$ as well.

Now, observe that the simple objects of $\PP'_{a,b}(1)$ are skyscraper sheaves $k(x)$ ($x\in\P^2$)
and objects of the form $i_*\GG[1]$, for $\GG\in\Coh \P^2$ a locally-free slope-stable sheaf on $\P^2$ with $\mu(\GG)=B$. Indeed, this can be proved in precisely the same way as Lemma \ref{lem:geom-stability}, (\ref{stableinP1}).
In particular, $\HH^0(A_j)$ is a torsion sheaf of dimension zero.

Let $B_j$ be the cokernel in $\Coh_0^{\sharp(B)}$ of $A_j\into F$; then we have an exact sequence
\[
0\to\FF_0 \stackrel{f_j}{\to}\HH^{-1}(B_j)\to \HH^0(A_j) \stackrel{g_j}{\to}\HH^0(F).
\]
The cokernel $\cok f_j$ is zero-dimensional; since $\FF_0$ is fixed 
$\HH^{-1}(B_j)$ is pure of dimension $2$, the length of $\cok f_j$ is bounded. As the length
of the image $\im g_j$ is also bounded, we get a bound on the length of $\HH^0(A_j)$.
At the same time, if $D_j$ denotes the cokernel in $\Coh_0^{\sharp(B)}$ of $A_j\into A_{j+1}$, we must have
$\HH^{-1}(D_j)=0$,  and thus an inclusion
$\HH^0(A_j)\into\HH^0(A_{j+1})$.
Hence, for $j\gg0$, $\HH^0(A_j)\cong\HH^0(A_{j+1})$ and so \eqref{eq:Bonn17} stabilizes.
\end{Prf}

Notice that, in the assumptions of Lemma \ref{lem:HNFiltrationsDiscrete}, the pair $(Z_{a,b},\Coh_0^{\sharp(B)})$ defines a locally-finite stability condition on $\DD_0$.
Indeed this follows immediately from \cite[Lemma 4.4]{Bridgeland:K3}.

\medskip

In the rest of this section, we will use Bridgeland's deformation result to
extend the existence of Harder-Narasimhan filtrations to the case
where $a, b$ are not rational. In order to make the deformation effective,
we need to bound the metric $\|\cdot\|_{\sigma_{a,b}}$ on $\Hom(K(\DD_0), \C) \cong \C^3$ defined in equation
\eqref{eq:def-metric} relative to an arbitrarily chosen metric, with the bound
depending continuously on $a, b \in G$.

We define the following functions $G \to \R_{\ge 0}$:
\begin{align*}
\gamma_{\min}(a, b) &:= \inf \stv{\abs{\gamma_{a,b}(t)}} {t \in \R} \\
E_{\min}(a,b) &:= \inf \stv{\abs{ \frac{Z_{a,b}(\EE)}{r(\EE)} + t}}
			{\text{$t \in \R_{\ge 0}$, $\EE$ exceptional v. bundle}} \\
S_{\min}(a,b) &:= \min( \gamma_{\min}(a,b), E_{\min}(a,b)).
\end{align*}

\begin{Lem}\label{lem:Smin}
The function $S_{\min}$ is continuous and satisfies
\[ 0 < S_{\min}(a,b) \le \inf \stv{\abs{ \frac{Z_{a,b}(\FF)}{r(\FF)} + t}}
			{t \in \R_{\ge 0},\, \FF \in \Coh_0\, \mathrm{slope}\text{-}\mathrm{stable}} \]
for all $(a, b) \in G$.
\end{Lem}
\begin{Prf}
The path $\gamma_{a,b}(t)$ depends continuously on $a, b, t$ and has the properties
$\gamma_{a,b}(t)\neq0$ for all $a, b, t$, and $\lim_{t \to \pm \infty} \abs{\gamma_{a,b}(t)} = +\infty$. It follows that $\gamma_{\min}$ is
a positive continuous function.

Since $Z_{a, b}(\EE) \not\in \R_{\leq 0}$, the term 
$\inf \stv{\abs{ \frac{Z_{a,b}(\EE)}{r(\EE)} + t}} {t \in \R_{\ge 0}}$ is positive for every
exceptional vector bundle $\EE$ on $\P^2$. Further, Theorem \ref{thm:DP} together with the
computation of Lemma \ref{lem:ineq} shows that every accumulation point of the set
\[ \stv {\frac{Z_{a,b}(\EE)}{r(\EE)}}
	{\text{$\EE$ exceptional v. bdle}} \]
is contained in $S_{a,b}$. Hence $E_{\min}$ is also a positive continuous function.

It remains to prove $S_{\min}(a, b) \le \abs{ \frac{Z_{a,b}(\FF)}{r(\FF)} + t}$ for all $\FF, t$. 
It holds by definition when $\FF$ is an exceptional vector bundle. Otherwise, the claim follows
as $\frac{Z_{a,b}(\FF)}{r(\FF)} + t$ is contained in $S_{a,b}$.
\end{Prf}

Let $\abs{\cdot}_\infty$ be the supremums-norm on $K(\DD_0)\otimes \R \cong \R^3$ in the coordinates
$(r, d, c)$, i.e., $\abs{(r,d,c)}_\infty=\max\{\abs{r},\abs{d},\abs{c}\}$.
Let $M(a, b)$ be the matrix
\[ M(a,b) := \begin{pmatrix} 1 & 0 & 0 \\ 
		\Im b & \Im a & 0 \\
		\Re b & \Re a & -1 
\end{pmatrix}, \]
and let $N(a, b) := \left\|M(a,b)^{-1}\right\|_\infty$ be the norm of
its inverse, where $\left\|\cdot\right\|_\infty$ is the operator norm with
respect to the supremums-norm on $\R^3$. We claim the following estimate,
which is the ``support property'' discussed in Proposition \ref{prop:SupportProperty}
with an explicit constant:

\begin{Lem}\label{lem:estimate}
If $E \in \Coh_0^{\sharp(B)}$ is $Z_{a,b}$-stable, then
\[
\frac{\abs{Z_{a,b}(E)}}{\abs{E}_\infty}
	\ge \frac{\min\left(S_{\min}(a, b), 1\right)}{N(a,b)}.
\]
\end{Lem}
We first show how to conclude the proof of Theorem \ref{thm:geom-stability} from the lemma:

\begin{Cor}
There exists a geometric stability condition $\sigma_{a,b}$ for arbitrary pairs
$(a,b) \in G \subset \C^2$.
\end{Cor}
\begin{Prf}
Let $V\subset U$ be the subset of geometric stability conditions $\sigma = (Z,\PP)$ such that all
skyscraper sheaves $k(x)$ are stable of phase $1$ with $Z(k(x)) = -1$. For any such stability
condition, the central charge is of the form $Z = Z_{a,b}$ of equation \eqref{eq:Zab}, and thus
$\ZZ$ induces a map $\ZZ_V \colon V \to \C^2, \sigma \mapsto (a,b)$.

As proved at the end of Section \ref{sec:Constraining}, the heart of such a stability condition is
uniquely determined; thus $\ZZ_V$ is injective, and so it is a
homeomorphism onto its image.  Using Corollary \ref{cor:openness} and the deformation property, we see
that the image of $\ZZ_V$ is open in $\C^2$. By Lemma \ref{lem:HNFiltrationsDiscrete}
and Lemma \ref{lem:estimate} (which shows that the stability conditions satisfy condition
\eqref{enum:full} in Definition \ref{def:GeomStability}), it contains
the dense subset of $(a,b) \in G$ such that $B = - \frac{\Im b}{\Im a} \in \Q$ is rational. Hence it
suffices to prove that the image is closed in $G$.

Assume the contrary, and that $(a, b) \in G$ are in the boundary of $\ZZ_V$.
By Lemma \ref{lem:Smin} the function $S_{\min}$ is continuous and positive.
It follows that for all $(a', b') \in G$ sufficiently
close to $(a,b)$, we have
\[ \left\|Z_{a,b} - Z_{a',b'}\right\|_{\infty}
 < \sin\left(\frac{\pi}8\right) \frac{\min\left(S_{\min}(a', b'), 1\right)}{N(a',b')}. \]
By the definition of $\|\cdot\|_{\sigma_{a', b'}}$ (see equation \eqref{eq:def-metric}) and
Lemma \ref{lem:estimate}, this implies
\[ \left\|Z_{a,b} - Z_{a',b'}\right\|_{\sigma_{a',b'}} < \sin\left(\frac{\pi}8\right) \]
for all such $(a', b')$ for which a geometric stability condition $\sigma_{a', b'}$ exists. 

By Bridgeland's effective deformation result (see Theorem \ref{thm:B-deform}) there exists a stability condition $\sigma_{a,b} = (Z_{a,b},\PP_{a,b})$ in the
neighborhood of $\sigma_{a', b'}$. By choosing $(a', b')$ appropriately, we may assume that
$\sigma_{a,b}$ is on one of the walls in the sense of Corollary \ref{cor:openness}; in particular,
$k(x)$ is semistable, and there is an inclusion $E \into k(x)$ in $\PP_{a,b}(1)$ with $E$ being
stable. In particular, $0 = \Im Z_{a,b}(E) = \Im a \cdot d(E) + \Im b \cdot r(E)$. Since
$k(x)$ is stable with respect to $\sigma_{a', b'}$, we have $\Im Z_{a', b'}(E) \neq 0$,
and thus we have $r(E) \neq 0$ or $d(E) \neq 0$. Since $\Im a > 0$, it follows that $r(E) \neq 0$;
but then $B = - \frac{\Im b}{\Im a} = \frac{d(E)}{r(E)} \in \Q$, and so we already know that
there exists a geometric stability condition $\sigma \in V$ with $\ZZ_V(\sigma) = (a,b)$.
\end{Prf}

\begin{Prf} (Lemma \ref{lem:estimate})
Writing $r = r(E)$ etc., we have
\begin{align*} \abs{E}_\infty &= \max(r,d,c) \\
& = \abs{M(a,b)^{-1} \cdot \bigl(r, \Im a \cdot d + \Im b \cdot r,
	    -c + \Re a \cdot d + \Re b \cdot r\bigr)}_\infty \\
& \le N(a,b) \cdot \abs{\left(r, \Im Z_{a,b}(E), \Re Z_{a,b}(E)\right)}_\infty \\
& \le N(a,b) \cdot \max(\abs{r}, \abs{Z_{a,b}(E)})
\end{align*}
Thus the claim follows if we can show 
\[ \abs{\overline Z(E)} > S_{\min(a, b)}
\] 
where $E \in \Coh_0^{\sharp(B)}$ is any $Z_{a,b}$-stable objects with non-zero rank, and where we
wrote $\overline Z(E) = \frac{Z_{a,b}(E)}{\abs{r(E)}}$.

Assume first that $r(H^0(E)) \ge r(H^{-1}(E))$, with $r(H^0(E)) > 0$. We have
\[
\Im Z_{a,b}(E) \ge \Im Z_{a,b}(H^0(E)),
\]
as $E \onto H^0(E)$ is a quotient in $\Coh_0^{\sharp(B)}$ and $Z_{a,b}$ is a stability function,
and we have $\abs{r(E)} \le r(H^0(E))$; together, they show
\[ \Im \overline Z(E) \ge  \Im \overline Z(H^0(E)). \]
Let $H^0(E) \onto \FF$ be a semistable quotient with $\mu(H^0(E)) \ge \mu(\FF) > B$ (such a quotient
always exists due to the existence of HN-filtrations for slope stability). Then
\begin{align*}
 \Im \overline Z(H^0(E))
& = \frac {\Im a \cdot d(H^0(E)) + \Im b \cdot r(H^0(E))} {r(H^0(E))} 
 = \Im a \left(\mu(H^0(E)) - B \right) \\
& \ge \Im a \left(\mu(\FF) - B\right) = \Im \overline Z(\FF)
\end{align*}
(where we used the assumption $\Im a > 0$ in the inequality).

On the other hand, as $E$ is $Z_{a,b}$-semistable and $\FF \in \Coh_0^{\sharp(B)}$, and so the phase
of $Z_{a,b}(\FF)$ is at least as big as the phase of $Z_{a,b}(E)$. Hence the line segment
connecting 0 and $\overline Z(E)$ intersects the ray $\overline Z(\FF) + t, t \ge 0$ (see Figure
\ref{fig:Z-proof}). By Lemma \ref{lem:Smin}, this implies the claim.

\begin{figure}
        \includegraphics{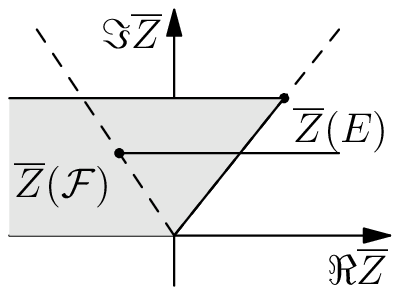}
\caption{Location of $\overline Z(\FF)$ relative to $\overline Z(E)$}
\label{fig:Z-proof}
\end{figure}

A dual argument holds in case $r(H^{-1}(E)) > r(H^0(E))$, by considering a
slope-stable sheaf $\FF \into H^{-1}(E)$ with $B \ge \mu(\FF) \ge \mu(H^{-1}(E))$. 
Finally, when $H^{-1}(E)$ is zero and $H^0(E)$ has rank zero, there is nothing to prove.
\end{Prf}

\begin{Rem}
The methods used in the last part of this section also apply in the situation
\cite[Section 2]{Aaron-Daniele}: the stability conditions constructed there for rational divisors
$D, F$ deform to produce stability condition for arbitrary $\R$-divisors $D, F$ with $F$ ample.
\end{Rem}

\section{Boundary of the geometric chamber}\label{sec:boundary}

In this section we will show that the set of boundary walls of
the geometric chamber $U$ can be described explicitly
using exceptional vector bundles $\EE$ on $\P^2$. For any such $\EE$, the
push-forward $i_*\EE$ is a spherical object in $\DD_0$. We denote
by ${\ST}_\EE \colon \DD_0 \to \DD_0$ the spherical twist at $i_* \EE$.
If $r$ is the rank of $\EE$ and $x \in \P^2$, we also write
$\EE^x$ for the kernel of the natural map
$i_*\EE^{\oplus r} \onto k(x)$.

The goal of this section is to prove the following theorem:

\begin{Thm}\label{thm:boundary}
For every exceptional vector bundle $\EE$ there exist two codimension one
walls $W_\EE^+, W_\EE^- \subset \partial U$ with the following properties:
\begin{enumerate}
\item 
Stability conditions in $W_\EE^+$ are characterized by the property that
$i_*\EE$ and all skyscraper sheaves $k(x)$ are semistable of the same phase
$\phi$, with $i_*\EE$ being a subobject of $k(x)$ in $\PP(\phi)$; at a
general point of $W_\EE^+$,
the Jordan-H\"older filtration of any skyscraper sheaf $k(x)$ is given by
\begin{equation}				\label{eq:kx-filtration1}
 i_*\EE^{\oplus r} \to k(x) \to \EE^x[1]. \end{equation}
Similarly, we have $\sigma \in W_\EE^-$ if $i_*\EE[2]$ is semistable of the
same phase $\phi$ as $k(x)$, and $i_*\EE[2]$ is a quotient of $k(x)$ in
$\PP(\phi)$. At a general point in $W_\EE^-$, the Jordan-H\"older
filtration of $k(x)$ is given by
\begin{equation}				\label{eq:kx-filtration2}
 {\ST}_\EE^{-1}(\EE^x[1]) \to k(x) \to i_*\EE^{\oplus r}[2]. \end{equation}
\item We have $W_\EE^+ = \overline{U} \cap {\ST}_\EE \left(\overline{U}\right)$,
i.e., $W_\EE^+$ is the wall between $U$ and ${\ST}_\EE(U)$; similarly,
$W_\EE^- = \overline{U} \cap {\ST}_\EE^{-1} \left(\overline{U}\right)$.
\end{enumerate}
There are no other walls in $\partial U$.
\end{Thm}

\begin{Cor}\label{cor:ConnectedComponent}
The translates of $\overline{U}$ under the group of autoequivalences generated by ${\ST}_\EE$ cover the
whole connected component $\Stab^{\dag}(\DD_0)$ of $U$ in the space of stability conditions.
\end{Cor}
\begin{Prf}
Let $\alpha \colon [0,1] \to \Stab(\DD_0)$ be a path of stability conditions with $\alpha(0) \in U$.
By Proposition \ref{prop:chambers}, there exists a finite set of walls $W^{[k(x)]}_i$ intersecting
$\alpha$, such that the set of stable objects of class $[k(x)]$ is constant in the complement of
the intersection points. We may also assume that $\alpha$ intersects each wall transversely and in a
generic point of the wall. Using the above theorem, it follows by induction that every open interval in
the complement is contained in the translate of $U$ under a sequence of spherical twists.
\end{Prf}

Most of the existing proofs of statements similar to the above claims are in the situation of a
Calabi-Yau 2-category. In that situation, Lemma 5.2 of \cite{Bridgeland:K3} applies, which
guarantees, via an Euler characteristic computation, that every non-trivial Harder-Narasimhan
filtration of a skyscraper sheaf $k(x)$ contains a spherical object. In a Calabi-Yau 3-category,
spherical objects cannot be characterized via their Euler characteristic among stable objects.
Instead, our proof is obtained by a direct geometric analysis of the boundary of $U$.

Consider a stability condition $(Z, \PP)$ in the boundary 
$\partial U$ of $U$. Since the skyscraper sheaves $k(x)$ are semistable,
$Z$ satisfies $Z(k(x)) \neq 0$ and, up to the action of $\C$, we
can still assume $k(x) \in \PP(1)$ and $Z = Z_{a,b}$ for
$a, b \in \C$ as in equation \eqref{eq:Zab}.
By Bridgeland's deformation result and Theorem 
\ref{thm:geom-stability}, 
$a,b \in \C$ must satisfy one of the following conditions:
\begin{description}
\item[Case $a$] $\Im a = 0$	
\item[Case $E$] 
$\Im a > 0$, there exists an exceptional vector bundle
of slope $B = - \frac {\Im b}{\Im a}$, and we have
$\delta_\infty^{DP}(B) > - \Re b  - B \cdot \Re a + \frac 12 B^2 > \Delta_B$.
\item[Case $\delta_\infty^{DP}$] 
$\Im a > 0$ and 
$- \Re b  - B \cdot \Re a + \frac 12 B^2 = \delta_\infty^{DP}(B)$.
\end{description}

We begin by showing that ``Case $\delta_\infty^{DP}$'' cannot exist.

\begin{Lem}\label{lem:constantheart}
Let $\sigma_t = (Z_t,\PP_t)$ for $t \in I \subset \R$ be a path
in the space of stability conditions such that $\Im Z_t$ is constant. 
Then $\PP_t((0,1])$ and $\PP_t(1)$ are constant, too.
\end{Lem}
\begin{Prf}
Let $t_1, t_2 \in I$ be such that 
$\sigma_{t_1}$ and $\sigma_{t_2}$ are close with respect
to the metric on $\Stab (\DD_0)$; to be specific,
we assume $d(\sigma_{t_1}, \sigma_{t_2}) < \frac 18$.
Given $\phi \in (0,1]$, the objects $E \in \PP_{t_2}(\phi)$ 
can be characterized as the $Z_{t_2}$-stable objects in the quasi-abelian 
category $\PP_{t_1}((\phi - \frac 18, \phi + \frac 18))$ (see
\cite[Section 7]{Bridgeland:Stab}). We want to show $E \in \PP_{t_1}((0,1])$.

In case $\frac 18 \le \phi \le \frac 78$ we are done. If
$\phi \in (\frac 78, 1]$, then for any $A \in \PP_{t_1}((1, \phi + \frac 18))$
we have $\Im Z_{t_2} (A) = \Im Z_{t_1}(A) < 0$, and thus the phase
of $A$ with respect to $Z_{t_2}$ is bigger than one. Hence $A$ cannot
be a subobject of $E$. By considering the Harder-Narasimhan filtration of
$E$ with respect to $\PP_{t_1}$, this implies that
$E \in \PP_{t_1}((\phi - \frac 18, 1])$.
A similar argument applies for $0 < \phi < \frac 18$.

It follows that $\PP_{t_2}((0,1]) \subset \PP_{t_1}((0,1])$, and thus they
must be equal. The claim about $\PP_t(1)$ follows easily.
\end{Prf}

\begin{Lem}
There are no stability conditions $\sigma = (Z_{a,b}, \PP)$ in $\partial U$
such that $\Im a > 0$ and
$-\Re b  - B \cdot \Re a + \frac 12 B^2 = \delta_\infty^{DP}(B)$.
\end{Lem}
\begin{Prf}
Due to Corollary \ref{cor:openness}, we may assume that
$B$ is irrational.

Consider the path $\sigma_t = (Z_{a(t), b(t)}, \PP_t), t \in [0,1]$ in
$\Stab (\DD_0)$ starting at
$\sigma$ induced by deforming $Z$ as $Z_{a(t), b(t)}$, with $a(t) = a$
constant, and $b(t) = b + \epsilon t$. Due to the form of the inequalities in Theorem \ref{thm:geom-stability}, the stability conditions $\sigma_t$ for $t> 0$
are geometric.

By the previous lemma, it follows that
$\PP_0((0,1]) = \PP_1((0,1]) = \Coh_0^{\sharp(B)}$ and, since $B$ is irrational, that
$\PP_0(1)$ is generated by the skyscraper sheaves $k(x)$.
The sheaves $k(x)$ have no subobjects in this category, thus they are stable.
By Corollary \ref{cor:openness}, this contradicts the assumption that
$\sigma$ is in the boundary of $U$.
\end{Prf}

The following lemma deals with ``Case $a$'':
\begin{Lem}
Let $(Z_{a,b}, \PP) \in \partial U$ be a stability condition with
$\Im a = 0$. Then $\Im b = 0$, i.e., the image of $Z$ is contained in the real
line.
\end{Lem}
\begin{Prf}
Writing inequality (\ref{ineq:genslope}) without denominators we get
\[ (\Im a)^2 \Re b > \Im a \cdot \Im b \cdot \Re a - (\Im a)^2 \delta_\infty^{DP}(B)
			+ \frac 12 (\Im b)^2 
\]
By continuity, this implies $(\Im b)^2 \le 0$ on the boundary with
$\Im a = 0$.
\end{Prf}

In particular, the part of the boundary with $\Im a = 0$ is a codimension two
subset; by Corollary \ref{cor:openness}, any such point is contained in the closure of the
$\Im a > 0$-part of a wall.

We will now consider boundary stability conditions in 
``Case $E$''.  Let $a, b, B$ be as in the assumption, and let $\EE$ be
the unique exceptional vector bundle on $\P^2$ of slope $B$ (see Theorem \ref{thm:DP-A}); then $Z_{a,b}(i_*\EE) \in (-1, 0)$ and
$Z_{a, b}(i_* \GG) \in \R_{> 0}$ for any other slope-stable sheaf $\GG$ on
$\P^2$ of slope $B$. This suggests that $i_*\EE$ is semistable,
of phase $\phi(i_*\EE) = \pm 1$ (depending on
whether $Z(i_*\EE)$ approaches the real line from above or below when
we approach $Z_{a,b}$ by geometric stability conditions); and in the case $\phi(i_*\EE) = +1$, the $t$-structure $\PP((0,1])$ should be given as in the following proposition.

We will prove this by constructing the stability conditions in the boundary directly, and prove that they deform to geometric stability conditions.

\begin{Prop}				\label{prop:E-tstruct}
Let $B$ be the slope of an exceptional vector bundle $\EE$ on $\P^2$. Then there
is a torsion pair $(\TT^\EE, \FF^\EE)$ on $\Coh_0$ where
\begin{itemize}
\item $\TT^\EE$ is the extension-closed subcategory of $\Coh_0$ generated by torsion
sheaves, by slope-semistable sheaves of slope $\mu > B$, and by $\EE$, and
\item $\FF^\EE$ is generated by slope-semistable sheaves $\GG$ of slope $\mu \le B$
that also satisfy $\Hom(i_* \EE, \GG) = 0$.
\end{itemize}
\end{Prop}
\begin{Prf}
From the construction, it is clear that $\Hom(\TT^\EE,\FF^\EE) = 0$.

Now given any $\GG \in \Coh_0$, let $\GG_{> B} \into \GG \onto \GG_{\le B}$
the unique short exact sequence with $\GG_{> B} \in \Coh_0^{> B}$ and
$\GG_{\le B} \in \Coh_0^{\le B}$.  Let $V = \Hom(i_*\EE, \GG_{\le B})$ and let
$\FF$ be the cokernel such that the following sequence is exact on the
right:
\[ V \otimes i_*\EE \stackrel{f}{\longrightarrow} \GG_{\le B} \to \FF \to 0 \]
We claim that $\FF \in \FF^\EE$, and that the kernel
$\TT$ of the composition $\GG \onto \GG_{\le B} \onto \FF$
lies in $\TT^\EE$.

The image of the evaluation map $f$ is
slope-semistable of slope $B$; hence so is the kernel of $f$. Since
$i_*\EE$ is stable, the kernel is of the form $i_* \EE \otimes V'$ for some
$V' \subset V$; by the definition of $V$, this forces $V' = 0$, i.e.,
the above sequence is exact. This shows that $\TT$ is an extension
of $\GG_{> B}$ and $i_*\EE$.

Since $\Hom(i_*\EE, i_*\EE) = \C$ and $\Ext^1(i_*\EE, i_*\EE) = 0$, the long
exact sequence associated to $\Hom(i_*\EE, \blank)$ shows
that $\Hom(i_*\EE, \FF) = 0$, and thus $\FF \in \FF^\EE$ as desired.
\end{Prf}

A similar result as Proposition \ref{prop:E-tstruct} is also in \cite[Prop.\ 2.7]{Yoshioka:StabilityFM}.

We continue to assume that $a, b \in \C$ satisfy $B = - \frac{\Im b}{\Im a}$ and the inequality of
boundary ``Case $E$''. Let $\Coh_0^\EE = \langle \TT^\EE, \FF^\EE[1] \rangle$ be the $t$-structure given
by tilting at the torsion pair of Proposition \ref{prop:E-tstruct}.  From the previous discussion,
it follows that $Z_{a, b}$ is a stability function for $\Coh_0^\EE$.

\begin{Prop}\label{prop:HNFiltrsBoundary}
Harder-Narasimhan filtrations exist for the stability function $Z_{a,b}$ on $\Coh_0^\EE$.
\end{Prop}
\begin{Prf}
The proof is similar to Lemma \ref{lem:HNFiltrationsDiscrete}. 
Indeed, since $B = -\frac{\Im b}{\Im a}$ is the slope of $\EE$, it is rational; hence the imaginary part
of $Z_{a,b}$ is discrete, and we can again apply Proposition \ref{prop:HNFiltrationsDiscrete}.

First note that for any $A \in \PP'_{a,b}(1)$, the sheaf $\HH^0(A)$ is in the category
extension-generated by $k(x), x \in \P^2$, and by $i_*\EE$: otherwise
$\HH^0(A) \in \TT^\EE$ would have $0 < \Im Z(\HH^0(A)) \le \Im Z(A)$.
To adapt the proof
of Lemma \ref{lem:HNFiltrationsDiscrete} (where instead we knew that $\HH^0(A)$ is a
zero-dimensional torsion sheaf), we replace all arguments using the length of a maximal
zero-dimensional subsheaf by using the function 
\[ e \colon \Coh_0 \to \Z, \quad e( \cdot) := \dim(\Hom(i_*\EE, \cdot)). \]
In many respects, it has the same formal properties needed (e.g.\ subadditivity on short exact sequences),
and the proof goes through:

Assume that $F\in\Coh_0^\EE$ has an infinite sequence
\begin{equation}\label{eq:Bonn26}
0=A_0\subset A_1\subset\ldots\subset A_j\subset A_{j+1}\subset\ldots\subset F,
\end{equation}
of subobjects with $A_j\in\PP'_{a,b}(1)$.
Denote by $B_j$ the cokernel in $\Coh_0^\EE$ of $A_j\into F$.

By arguing as in Lemma \ref{lem:HNFiltrationsDiscrete}, we can assume $\HH^{-1}(A_j)=\FF_{-1}$ and $\HH^0(B_j)=\GG_1$, for some $\FF_{-1},\GG_1\in\Coh_0$ and for all $j$.
Let $\FF_0$ be the cokernel in $\Coh_0$ of $\FF_{-1}\into\HH^{-1}(F)$ and let $\GG_0$ be the kernel of $\HH^0(F)\onto\GG_1$.
Then we have an exact sequence
\[
0\to\FF_0\stackrel{g}{\to}\HH^{-1}(B_j)\to\HH^0(A_j)\stackrel{f}{\to}\GG_0\to0, 
\]
and we let $\KK_j = \ker f = \cok g$. By the long exact $\Hom$-sequences, we have
\[
e(\HH^0(A_j))  \le e(\KK_j) + e(\GG_0)
\le e(\HH^{-1}(B_j)) + \dim \Ext^1(i_*\EE,\FF_0) + e(\GG_0).
\]
By definition of $\FF^\EE$, we have $\Hom(i_* \EE, \HH^{-1}(B_j)) = 0$, and thus 
$e(\HH^0(A_j))$ is bounded.

Now consider a filtration step of \eqref{eq:Bonn26}, and 
let $D_j$ be the cokernel in $\Coh_0^\EE$ of $A_j\into A_{j+1}$.
From the exact sequence
\[
0\to\HH^{-1}(D_j)\to\HH^0(A_j)\to\HH^0(A_{j+1})\stackrel{\phi}{\to}\HH^0(D_j)\to0,
\]
we argue similarly as before:
\begin{align*}
 e(\ker \phi) = e(\HH^{-1}(D_j)) + e(\ker \phi) \ge e(\HH^0(A_j))
\end{align*}
On the other hand, as $\TT^\EE$ is closed under quotients, we have
$\ker \phi \in \TT^\EE$, and thus $\ker \phi \in \PP'_{a,b}(1)$ and
$\Ext^1(i_*\EE, \ker \phi) = 0$. Thus we have
\[
e(\HH^0(A_{j+1})) = e(\ker \phi) + e(\HH^0(D_j)) \ge e(\ker \phi) \ge e(\HH^0(A_j))
\]
with equality only if $\HH^0(D_j) = 0$. By the boundedness established above,
we do have $\HH^0(D_j) = 0$ for $j \gg 0$; but then
$\HH^0(A_j) \onto \HH^0(A_{j+1})$ stabilizes as $\Coh_0$ is Noetherian.
\end{Prf}

We denote by $\WW_\EE^+$ the set of stability conditions constructed in the above proposition.
The same methods as in Section \ref{sec:constructing} show that they are full, i.e., that they
satisfy condition \ref{enum:full} of Definition \ref{def:GeomStability}; alternatively, this follows
from Corollary \ref{cor:BoundaryIsAlgebraic}, where it is shown that they can also be constructed
as a stability condition whose heart of finite length are representations of a finite quiver.

\begin{Prop}
For any $\sigma \in \WW_\EE^+$, the skyscraper sheaves $k(x)$ are $\sigma$-semistable
with Jordan-H\"older filtration given as in \eqref{eq:kx-filtration1}.
Their images under the spherical twist at
$i_*\EE$ are also $\sigma$-semistable, of the same phase as $k(x)$, with Jordan-H\"older
filtration given by 
\begin{equation}				\label{eq:TEkx}
 \EE^x[1] \to {\ST}_\EE(k(x)) \to i_*\EE^{\oplus r}. \end{equation}
\end{Prop}

\begin{Prf}
The sheaves $\EE^x$ are slope-semistable of the same slope as $\EE$, and the
long exact $\Hom$-sequence shows $\Hom(i_*\EE, \EE^x) = 0$. Hence
$\EE^x \in \FF^\EE$, and $\EE^x[1] \in \PP(1)$, and we indeed have a
short exact sequence as in (\ref{eq:kx-filtration1}) in $\PP(1)$.
We claim that $\EE^x[1]$ is $Z_{a,b}$-stable, i.e., that there are no
non-trivial short exact sequences $M \into \EE^x[1] \onto N$ in the abelian category
$\PP(1)$:

Let $C$ denote the kernel of the composition $k(x) \onto \EE^x[1] \onto N$ in $\PP(1)$.
By the long exact cohomology sequence, $C$ is isomorphic to a sheaf $\CC$. As observed in the proof
of Proposition \ref{prop:HNFiltrsBoundary}, $\CC = \HH^0(C)$ lies in the category extension-generated by
$i_*\EE$ and skyscraper sheaves $k(x)$ for $x \in \P^2$.
We claim that due to $\Ext^1(k(x), i_*\EE) = \Ext^1(i_*\EE, i_*\EE) = 0$, this implies more strongly
that there is a short exact sequence
$\TT \into \CC \onto i_*\EE^{\oplus k}$ for some zero-dimensional torsion sheaf $\TT$.
Indeed, by induction on the length of the Jorder-H\"older filtration of $\CC$ in this finite-length
category, we may assume that there is a sheaf $\CC' \in \PP(1)$ with
a short exact sequence $\TT' \into \CC' \onto i_*\EE^{\oplus k'}$, such that $\CC$ is an extension
of $\CC'$ by a simple object, i.e., there is a short exact sequence either of the form
$k(x) \into \CC \onto \CC'$ or $i_*\EE \into \CC \onto \CC'$.
In the former case, the claim for $\CC$ follows immediately by considering the composition $\CC
\onto \CC' \onto i_*\EE^{\oplus k'}$; in the latter case, the vanishing of
$\Ext^1(\TT, i_*\EE) = 0$ implies that there is factorization $\TT' \into \CC$; due to
$\Ext^1(i_*\EE, i_*\EE) = 0$, the kernel is of the form $i_*\EE^{\oplus k'+1}$.

The composition $\TT \into \CC \into k(x)$ can only be injective
in $\PP(1)$ if it is injective as a map of sheaves; hence either $\TT \cong \CC \cong k(x)$,
or $\TT = 0$ and $\CC \cong i_*\EE^{\oplus k}$.

In the former case we have $N = 0$. In the latter case, note that the inclusion
$i_*\EE^{\oplus r} \into k(x)$ factors via 
$i_*\EE^{\oplus r} \into \CC \cong i_*\EE^{\oplus k} \into k(x)$ and induces an isomorphism
$\Hom(i_*\EE, i_*\EE^{\oplus r}) \cong \Hom(i_*\EE, k(x)) \cong \C^r$. Thus we must have $k= r$, and 
$N \cong \EE^x[1]$. So in both cases the exact sequence is trivial.

By using adjunction one sees that
$\RHom(i_*\EE, k(x)) = \C^r \oplus \C^r[-1]$. The long exact cohomology sequence
shows $\HH^0({\ST}_\EE(k(x))) \cong i_*\EE^{\oplus r}$ and
$\HH^{-1}({\ST}_\EE(k(x))) = \EE^x$, and so there is an exact triangle as in
(\ref{eq:TEkx}). This shows ${\ST}_\EE(k(x)) \in \PP(1)$, and that
(\ref{eq:TEkx}) is a Jordan-H\"older filtration.
\end{Prf}

We can deform $Z$ such that $Z(k(x)) = -1$ remains constant, $Z(i_*\EE)$ moves to the upper
half-plane, and $Z(\EE^x[1])$ moving to the lower-half plane; then by Lemma \ref{lem:semitostable},
all $k(x)$ become stable.  It follows that the closure $W_\EE^+$ of the orbit
of $\WW_\EE^+$ under the
action of $\C$ is a wall of $\partial U$. The objects $k(x)$ and ${\ST}_\EE(k(x))$ become stable on
opposite sides of the wall, and thus $W_\EE^+ \subset \overline{U} \cap {\ST}_\EE(\overline{U})$. If we apply
${\ST}_\EE^{-1}$ to $\WW_\EE^+$, we obtain a wall where the Jordan-H\"older filtration of $k(x)$ is given by
the image of (\ref{eq:TEkx}) under ${\ST}_\EE^{-1}$, which is indeed the exact triangle (\ref{eq:kx-filtration2}).

This finishes the proof of Theorem \ref{thm:boundary}.
Note that the proof also implies that two
such walls can only intersect at points where the image of the central charge is contained in a line
$e^{i\pi\phi}\cdot \R \subset \C$. In that case, the heart $\PP(\phi)$ of the associated t-structure
has finite length; in fact, it is one of the ``quivery'' stability conditions constructed in 
\cite{Bridgeland:stab-CY}, and used in the following section.

\begin{Lem}\label{lem:semitostable}
Let $E \in \DD_0$ and $\sigma \in \Stab^\dag(\DD_0)$ be a stability condition such that 
$E$ is $\sigma$-semistable, and assume that there is a Jordan-H\"older filtration
$M^{\oplus r} \into E \onto N$ of $E$ such that $M, N$ are $\sigma$-stable, $\Hom(E, M) = 0$, and
$[E]$ and $[M]$ are linearly independent classes in $K(\DD_0)$.
Then $\sigma$ is in the closure of the set of stability conditions where $E$ is
stable.
\end{Lem}
\begin{Prf}
Let $\phi$ be the phase of $M, E, N$ with respect to $\sigma = (Z, \PP)$.
By similar arguments as in the proof of Proposition \ref{prop:chambers}, we can show that for any stability condition $\sigma'$ sufficiently close to $\sigma$, $E$ can
only be destabilized by subobjects $F \into E$ in $\PP(\phi)$. Now let $\sigma' = (Z', \PP')$
be such a stability condition close by with $M, N$ stable and $\phi'(M) < \phi'(E) < \phi'(N)$.
Assume that $F \in \PP(\phi)$ is a stable destabilizing subobject of $E$ with respect to
$Z'$. If the image of the composition $F \to N$ is zero, then $F$ factors via
$M$, hence $\phi'(F) < \phi'(M)$. So $F \to N$ must be surjective; its kernel $G \in \PP(\phi)$ is a
subobject of $M^{\oplus r}$, and thus of the form $M^{\oplus k}$ for some $k < r$. Hence the quotient of
$F \into E$ is isomorphic to $M^{\oplus r-k}$, in contradiction to $\Hom(E, M) = 0$.
\end{Prf}

\section{Algebraic stability conditions}\label{sec:AlgebraicStability1}

In this section we study the open subset $\Stab_a$, introduced by Bridgeland in \cite{Bridgeland:stab-CY}, consisting of \emph{algebraic stability conditions}.
We first introduce open subsets $\Theta_\SSS$, associated to a collection of spherical objects $\SSS$, and study their boundary in $\Stab^\dag(\DD_0)$.
The subset of algebraic stability conditions will then be the union of all $\Theta_\SSS$.
Then we study in detail the relation between $\Stab_a$ and $U$. In particular, we show that $\Stab_a$ contains the boundary of $U$ (described in the previous section) and, vice versa, that the intersection of $U$ with $\Stab_a$ is strictly contained in $U$. In the next section we will apply all of this to prove Theorem \ref{thmi:sc}.

Let $\EEE=\{\EE_0,\EE_1,\EE_2\}$ be an exceptional collection of vector bundles on $\P^2$.
Recall that (see \cite{Goro-Ruda:Exceptional,Bondal}) a collection $\EEE$ of exceptional vector bundles is called \emph{exceptional} if $\Ext^p(\EE_j,\EE_i)=0$, for all $p$ and all $j>i$.
On $\P^2$, all exceptional collections also satisfy the vanishing (strong exceptional collection) $\Ext^p(\EE_i,\EE_j)=0$, for all $p>0$ and all $i<j$.
Moreover, every exceptional vector bundle on $\P^2$ is part of an exceptional collection of vector bundles.

The subcategory of $\DD_0$ generated by extensions by $i_*\EE_0[2]$, $i_*\EE_1[1]$, and $i_*\EE_2$
\[
\AA_{\EEE}:=\langle i_*\EE_2,i_*\EE_1[1],i_*\EE_0[2]\rangle
\]
is the heart of a bounded $t$-structure on $\DD_0$. By \cite{Bridgeland:TStruct}, the category $\AA_\EEE$ can also be described as the category of nilpotent modules over a certain algebra.

\begin{Def}\label{def:quivery}
A heart of a bounded $t$-structure on $\DD_0$ is called \emph{quivery} if it is of the form $\Phi(\AA_\EEE)$, for some exceptional collection $\EEE$ of vector bundles on $\P^2$ and for some autoequivalence $\Phi$ of $\DD_0$ given by composition of spherical twists associated to exceptional vector bundles.
\end{Def}

Notice that a quivery subcategory is of finite length, with simple objects $\Phi(i_*\EE_2)$, $\Phi(i_*\EE_1[1])$, and $\Phi(i_*\EE_0[2])$, which are also spherical in $\DD_0$. A quivery subcategory is called \emph{ordered} if it comes with an ordering of $S_0,S_1,S_2$ of its simple objects compatible with the requirement that $\Hom^k(S_j,S_l)=0$ unless $0\leq k\leq3$ and $j-l\equiv k$(mod $3$). A collection $\SSS=\{S_0,S_1,S_2\}$ of spherical objects of $\DD_0$ is called an \emph{ordered quivery collection} if it arises as an ordered collection of simple objects in an ordered quivery subcategory (which we will denote by $\AA_{\SSS}$).

By \cite[Theorem 4.11]{Bridgeland:TStruct}, we can define an action on the set of quivery ordered subcategories of $\DD_0$ of the affine braid group $B_3$, i.e., the group generated by elements $\tau_j$ ($j\in\Z_3$) and $r$ subject to the relations
\[
r\tau_jr^{-1}=\tau_{j+1},\qquad\tau_j\tau_{j+1}\tau_j=\tau_{j+1}\tau_j\tau_{j+1},\qquad r^3=1.
\]
Indeed, to define such an action is sufficient to set how the generators of $B_3$ act on the simple objects of a ordered quivery category:
\begin{align*}
\tau_1\{S_0,S_1,S_2\}&:=\{S_1[-1],{\ST}_{S_1}(S_0),S_2\}\\
r\{S_0,S_1,S_2\}&:=\{S_2,S_0,S_1\}.
\end{align*}
By \cite[Prop.\ 4.10]{Bridgeland:TStruct}, the image via $\tau_j$ of a quivery category is quivery as well, and thus the action is well-defined.
Notice, in particular, that
\[
\tau_2\{S_0,S_1,S_2\}=\{S_0,S_2[-1],{\ST}_{S_2}(S_1)\}.
\]

\begin{Rem}\label{rmk:GoroRuda}
Let $\EEE=\{\EE_0,\EE_1,\EE_2\}$ be an exceptional collection of vector bundles on $\P^2$.
Then
\[
\tau_1\{i_*\EE_2,i_*\EE_1[1],i_*\EE_0[2]\}=\{i_*\FF_2,i_*\FF_1[1],i_*\FF_0[2]\},
\]
where $\FFF=\{\FF_0,\FF_1,\FF_2\}$ is another exceptional collection of vector bundles on $\P^2$, called the \emph{left mutation} of $\EEE$ at $\EE_1$ (see \cite{Goro-Ruda:Exceptional,Bondal}). Similarly, for $\tau_2$ we have the \emph{left mutation} of $\EEE$ at $\EE_0$, for $\tau_1^{-1}$ we have the \emph{right mutation} of $\EEE$ at $\EE_2$, and for $\tau_2^{-1}$ we have the \emph{right mutation} of $\EEE$ at $\EE_1$.
Since all exceptional collections of vector bundles on $\P^2$ can be obtained by a sequence of mutations from $\EEE_1:=\{\OO_{\P^2}(-1),\Omega_{\P^2}(1),\OO_{\P^2}\}$, all ordered quivery subcategories can be obtained from
\[
\AA_1:=\AA_{\EEE_1}=\langle i_*\OO_{\P^2},i_*\Omega_{\P^2}(1)[1],i_*\OO_{\P^2}(-1)[2]\rangle
\]
by the action of $B_3$.
\end{Rem}

\begin{Def}\label{def:algebraicstability}
A stability condition $\sigma$ on $\DD_0$ is called \emph{algebraic} if there exists $M\in\widetilde{\GL}_2(\R)$ such that the heart of $\sigma\cdot M$ is quivery. Denote by $\Stab_a$ the subset of $\Stab(\DD_0)$ consisting of algebraic stability conditions.
\end{Def}

Using \cite[Prop.\ 4.10]{Bridgeland:TStruct} and \cite[Cor.\ 3.20]{Macri:Curves} it follows that $\Stab_a$ is an open connected $3$-dimensional submanifold of $\Stab(\DD_0)$.
Moreover, it is easy to construct stability conditions in $\Stab_a$ in which the skyscraper sheaves are all stable (for example, a stability condition with heart $\AA_1$ in which $\phi(i_*\OO_{\P^2})<\phi(i_*\Omega_{\P^2}(1)[1])<\phi(i_*\OO_{\P^2}(-1)[2])$). Hence $\Stab_a\subset\Stab^\dag(\DD_0)$, but the inclusion is strict (this can be deduced from Proposition \ref{prop:geomvsalg} and Remark \ref{rmk:Bonn1107}). Finally, by its own definition, $\Stab_a$ is invariant under the subgroup of the autoequivalences of $\DD_0$ which is generated by spherical twists ${\ST}_\FF$, with $\FF$ an exceptional bundle on $\P^2$.

\begin{Def}\label{def:theta}
Let $\SSS$ be an ordered quivery collection.
We denote by $\Theta_{\SSS}$ the open subset of $\Stab_a$ consisting of stability conditions whose heart is, up to the action of $\widetilde{\GL}_2(\R)$, equivalent to $\AA_\SSS$. With a slight abuse of notation, when $\EEE$ is an exceptional collection of vector bundles on $\P^2$, we denote by $\Theta_{\EEE}$ the open subset of $\Stab_a$ consisting of stability conditions whose heart is, up to the action of $\widetilde{\GL}_2(\R)$, equivalent to $\AA_\EEE$.
\end{Def}

\begin{Lem}\label{lem:thetaP2} The region
$\Theta_{\SSS} \subset \Stab(\DD_0)$ is characterized as the subset where $S_0, S_1, S_2$ are stable, and where their
phases $\phi_j := \phi(S_j)$ satisfy 
\begin{equation} \label{eq:phasesclose}
\abs{\phi_j - \phi_{j+1}} < 1 \quad \text{for $j= 0, 1, 2$}.\end{equation}
It is homeomorphic to
\[
\CC_{\SSS} = \left\{(m_0,m_1,m_2,\phi_0,\phi_1,\phi_2)\in\R^6\colon \text{\rm $m_j>0$ and \eqref{eq:phasesclose} holds for
all $j$} \right\}
\]
\end{Lem}
\begin{Prf}
Evidently $S_j$ are stable in $\Theta_\SSS$, and satisfy equation \eqref{eq:phasesclose}. Conversely, if $S_j$ are stable in $(Z, \PP)$
satisfying equation \eqref{eq:phasesclose}, then for $\phi$ slightly smaller than $\min \phi_j$ we have $\AA_\SSS \subset \PP((\phi, \phi+1])$, thus
$\AA_\SSS = \PP((\phi, \phi+1])$ and $(Z, \PP) \in \Theta_\SSS$.
\end{Prf}

Notice that, for later use, ${\ST}_{S_{j+1}}(S_{j})$ is an extension of $S_{j+1}$ by $m$ copies of $S_{j}$, where $m=\dim\Hom^1(S_{j+1},S_j)$.
Hence, its class in $K(\DD_0)$ is given by
\[
[{\ST}_{S_{j+1}}(S_j)]=m[S_j]+[S_{j+1}].
\]
Moreover, it belongs to $\AA_\SSS$ and, if $\phi(S_{j+1})>\phi(S_j)$ then it is also $\sigma$-stable.
A similar observation holds true for ${\ST}^{-1}_{S_{j}}(S_{j+1})$.

The next proposition generalizes \cite[Theorem 1.1]{Bridgeland:stab-CY}.

\begin{Prop}\label{prop:BoundaryTheta}
Let $\SSS=\{S_0,S_1,S_2\}$ be an ordered quivery collection.
Then the closure $\overline{\Theta}_{\SSS}$ of $\Theta_{\SSS}$ in $\Stab^\dag(\DD_0)$ is contained in $\Stab_a$.
\end{Prop}
\begin{Prf}
Let $\overline{\sigma}=(\overline{Z},\overline{\PP})\in\overline{\Theta}_{\SSS}\setminus\Theta_{\SSS}$ be the limit of a sequence $\{\sigma_s\}_{s\in\N}$, with $\sigma_s\in\Theta_{\SSS}$.
Then $S_0$, $S_1$, and $S_2$ are $\overline{\sigma}$-semistable; up to the action of $\C$, we have the following possibilities for their phases:
\begin{enumerate}
\item \label{enum:Zinline} The image of $\overline{Z}$ is a line in the plane.
\item \label{enum:jj+1} The image of $\overline{Z}$ is not a line in the plane and there exists $j$ such that $0=\overline{\phi}(S_j)=\overline{\phi}(S_{j+1})-1$ (here and in the sequel all the indices are taken modulo $3$).
\item \label{enum:jj-1} The image of $\overline{Z}$ is not a line in the plane and there exists $j$ such that $0=\overline{\phi}(S_j)=\overline{\phi}(S_{j-1})-1$.
\end{enumerate}

We begin with case (\ref{enum:jj+1}).
First of all notice that $0<\overline{\phi}(S_{j-1})<1$.
Moreover, up to the action of $\C$, we can assume every $\sigma_s$ has heart $\AA_{\SSS}$.

Let $P_m$ be the Kronecker quiver with $m = \dim\Hom^1(S_{j+1},S_j)>0$ , i.e., the quiver with two vertices and $m$ arrows from the first to the second vertex.
Consider the faithful functor $I \colon \Db(P_m)\to\mathrm{Tr}(S_j,S_{j+1})\subset\DD_0$, which maps the two simple quiver representations of $P_m$ corresponding to the two vertices respectively to $S_{j+1}$ and $S_j$. Here $\mathrm{Tr}(S_j,S_{j+1})$ denotes the triangulated subcategory of $\DD_0$ generated by $S_j$ and $S_{j+1}$.
For $s\gg 0$, the stability condition $\sigma_s$ induces a stability condition on $\mathrm{Tr}(S_j,S_{j+1})$, whose heart is the abelian category generated by extensions by $S_j$ and $S_{j+1}$. Now the functor $I$ restricted to mod-$P_m$ is full and faithful.
By \cite[Prop. 2.12]{MMS:inducing}, $\sigma_s$ induces a stability condition $I^{-1}\sigma_s$ in $\Db(P_m)$.
Hence, by \cite[Lemma 2.9]{MMS:inducing}, $I^{-1}\overline{\sigma}\in\Stab(P_m)$.
By \cite[Lemma 4.2]{Macri:Curves}, there exists an integer $k\in\Z$ such that $\overline{\sigma}\in\Theta_{\tau_{j+1}^k\SSS}$.
More explicitly, if $\tau_{j+1}^k\SSS=\{R_0,R_1,R_2\}$, then what we proved is that $R_j$ are stable with respect to $\overline{\sigma}$, and that we have Jordan-H\"older filtrations given by $S_{j-1}=R_{j-1}$ and
\begin{equation}\label{eq:Bonn1007}
\begin{split}
&R_{j}^{\oplus u_j}\to S_j[\epsilon]\to R_{j+1}^{\oplus v_j}\\
&R_{j}^{\oplus u_{j+1}}\to S_{j+1}[\epsilon-1]\to R_{j+1}^{\oplus v_{j+1}},
\end{split}
\end{equation}
where $\epsilon=0,1$ according to $k$.

For possibility (\ref{enum:jj-1}), we have similarly $0<\overline{\phi}(S_{j+1})<1$.
Then $\overline{\sigma}\in\Theta_{\tau_{j-1}\SSS}$ and $S_0$, $S_1$, and $S_2$ remain stable in $\overline{\sigma}$.

Finally, if the image of $Z$ lies in a line (case (\ref{enum:Zinline})), then we can deform $\overline{\sigma}$ in $\overline{\Theta}_{\SSS}$ in such a way to reduce to the situation of case (\ref{enum:jj+1}).
We can apply the previous procedure and find $g_1\in B_3$ such that $\overline{\sigma}\in\overline{\Theta}_{g_1\SSS}$.
If $\overline{\sigma}\in\Theta_{g_1\SSS}$ we have finished the proof. Assume not.
Then we continue and again deform $\overline{\sigma}$ in $\overline{\Theta}_{g_1\SSS}$ to reduce again to case (\ref{enum:jj+1}). We produce a new element $g_2\in B_3$ and so on.
This procedure must eventually terminate at a step $N$: indeed at every step, by \eqref{eq:Bonn1007}, we are constructing a filtration of $S_0$, $S_1$, and $S_2$ into $\overline{\sigma}$-semistable objects of the same phase. But $\overline{\sigma}$ is locally-finite. Hence at a certain point we produce a stable factor and so $\overline{\sigma}\in\Theta_{g_N\SSS}$, as wanted.
\end{Prf}

We can now study the relation of $\Stab_a$ with $U$.

\begin{Lem}\label{lem:StabilitySkyscrapers}
Let $\EEE$ be an exceptional collection of vector bundles on $\P^2$.
Then $\Theta_{\EEE}\cap U\neq\emptyset$ and it is connected.
\end{Lem}
\begin{Prf}
First of all notice that all skyscraper shaves $k(x)$ for $x \in \P^2$ belong to $\AA_\EEE$.
Consider the stability condition $\overline{\sigma}\in\Theta_\EEE$ with heart $\AA_{\EEE}$, whose simple objects have phases
\[
\overline{\phi}(i_*\EE_2)=\overline{\phi}(i_*\EE_1[1])=\overline{\phi}(i_*\EE_0[2])=1.
\]
Let $\EE_2^x$ be the kernel of $\EE_2^{\oplus r} \onto k(x)$ as in Section \ref{sec:boundary}; then 
$\EE_2^x[1] \in \AA_\EEE$, and in fact it is contained in the abelian category generated by extensions by $i_*\EE_1[1]$ and $i_*\EE_0[2]$.
We can deform $\overline{\sigma}$ slightly to a stability condition $\sigma\in\Theta_\EEE$ with
\begin{itemize}
\item $\phi(i_*\EE_2)=\phi(k(x))=\phi(\EE_2^x[1])$,
\item $\phi(i_*\EE_1[1])<\phi(i_*\EE_0[2])$.
\end{itemize}

We claim that $\EE_2^x$ is $\sigma$-stable.
Indeed, as in the proof of Proposition \ref{prop:BoundaryTheta}, we can consider the faithful functor $I\colon \Db(P_m)\to\mathrm{Tr}(i_*\EE_0,i_*\EE_1)$, where $m = \dim\Hom(i_*\EE_0,i_*\EE_1)>0$.
Then $\EE_2^x\cong I(\widetilde{\EE}_2^x)$, and to prove that $\EE_2^x$ is $\sigma$-stable is
equivalent to prove that $\widetilde{\EE}_2^x$ is $I^{-1}\sigma$-stable in $\Db(P_m)$. But the
stability of $\widetilde{\EE}_2^x$ follows immediately from \cite[Proposition 4.4]{King:QuiverStability}.

Hence $k(x)$ is $\sigma$-semistable, and its two Jordan--H\"older factors are $\EE_2^x[1]$ and $i_*\EE_2^{\oplus r_2}$, where $r_2$ is the rank of $\EE_2$.
By Lemma \ref{lem:semitostable}, $\sigma\in\Theta_\EEE\cap\overline{U}$, and so $\Theta_{\EEE}\cap U\neq\emptyset$ since $\Theta_\EEE$ is open.

To prove connectedness, we may first use the action by $\C$ to fix the phase of $k(x)$ to be 1 with $Z(k(x)) = -1$. Then
every class of a subobject of $k(x)$ gives a linear inequality for the imaginary part of $Z$, and thus $\Theta_\EEE \cap
U$ is cut out by a finite number of half-spaces.
\end{Prf}

\begin{Cor}\label{cor:BoundaryIsAlgebraic}
We have
\[
{\Stab}^{\dag}(\DD_0)={\Stab}_a\cup\bigcup\Phi(U),
\]
where the union is taken over all autoequivalences $\Phi$ of $\DD_0$ which belongs to the subgroup generated by spherical twists ${\ST}_\FF$ ($\FF$ an exceptional vector bundle on $\P^2$).
\end{Cor}
\begin{Prf}
By Corollary \ref{cor:ConnectedComponent}, we know that $\Stab^\dag(\DD_0)=\bigcup\Phi(\overline{U})$.
We only need to show that the boundary $\partial U$ is contained in $\Stab_a$.

Complete $\EE$ to an exceptional collection $\EEE=\{\EE_0,\EE_1,\EE_2=\EE\}$ of vector bundles on $\P^2$.
By the proof of Lemma \ref{lem:StabilitySkyscrapers}, there exists a stability condition
$\sigma = (Z, \PP) \in \Theta_{\EEE}\cap \WW_\EE^+$, where $\WW_\EE^+$ consists of the stability
conditions constructed by Propositions \ref{prop:E-tstruct} and \ref{prop:HNFiltrsBoundary}.
For such $\sigma$, the objects $i_*\EE_0[1]$ and $i_*\EE_1[1]$ belong to $\Coh_0^\EE$; combined
with Lemma \ref{lem:thetaP2} we get
\begin{equation}\label{eq:IneqPhases}
0 < \phi(i_*\EE_0[1]) < \phi(i_*\EE_1[1]) < 1 = \phi(i_* \EE_2).
\end{equation}
By definition $\WW_{\EE}^+$ is connected.
Moreover, $\Theta_{\EEE}\cap\WW_{\EE}^+\neq\emptyset$ is open in $\WW_{\EE}^+$.
We want to show it is also closed.

Let $\overline{\sigma}$ be a stability condition on the boundary
$\WW_{\EE}^+ \cap \partial \Theta_{\EEE}$.
Due to the inequalities \eqref{eq:IneqPhases}, the only inequality of Lemma \ref{lem:thetaP2} that
can become an equality for $\overline{\sigma}$ is 
\[
\overline{\phi}(i_*\EE_0[1])= \overline{\phi}(i_*\EE_1[1]).
\]
On the other hand, we have $\EE_2^x[1] \in \overline{\PP}(1)$ by construction of $\WW_{\EE}^+$ and, at the same
time, it lies in the abelian subcategory generated by $i_*\EE_0[2]$ and $i_*\EE_1[1]$. Therefore,
the central charge $\overline{Z}$ is contained in a line, i.e., $\overline{\sigma}\notin\WW_{\EE}^+$.

Hence $W_\EE^+ = \overline{\WW_{\EE}^+} \subset\overline{\Theta}_{\EEE}$ and, by Proposition \ref{prop:BoundaryTheta},
$W_\EE^+ \subset \Stab_a$.
For the case of the boundary of type $W_\EE^-$, simply observe that $W_\EE^-={\ST}_\EE^{-1}(W_\EE^+)\subset\Stab_a$.
By Theorem \ref{thm:boundary}, $\partial U\subset\Stab_a$, as wanted.
\end{Prf}

Notice that, in the proof of Corollary \ref{cor:BoundaryIsAlgebraic}, we actually showed that
\begin{equation}\label{eq:Utah29909}
\partial U\subset\bigcup\overline{\Theta}_\EEE,
\end{equation}
where the union is taken over all exceptional collections of vector bundles on $\P^2$. It follows that:
\begin{Rem} \label{rem:degenerate}
There is a one-to-one correspondence between quivery subcategories and loci in $\Stab^\dag(\DD_0)$ of codimension 2 where the image of the central charge is contained in a line.
\end{Rem}
Indeed, such a degenerate stability condition must lie, up to translation by spherical twists, in the boundary
$\partial U$. It has a unique heart (up to shifts), which must be $\AA_\EEE$ for some exceptional collection
$\EEE$.

\begin{Cor}\label{cor:IntersectionIsConnected}
$\Stab_a\cap U$ is connected.
\end{Cor}
\begin{Prf}
Let $\EEE=\{\EE_0,\EE_1,\EE_2\}$ be an exceptional collection of vector bundles on $\P^2$.
By Lemma \ref{lem:StabilitySkyscrapers}, $\Theta_\EEE\cap U$ is nonempty and connected.
We first claim that we can connect in $\Stab_a\cap U$ any stability condition in $\Theta_\EEE\cap U$ to a stability condition in $\Theta_{\EEE_1}\cap U$.

We proceed by induction on the length of a mutation from $\EEE_1$ to $\EEE$.
By Remark \ref{rmk:GoroRuda}, we need to show that a stability condition in $\Theta_{\tau_j^{\pm 1}\EEE}\cap U$, for $j=1,2$, can be connected to a stability condition in $\Theta_\EEE\cap U$.
Let $\sigma\in\Theta_{\EEE}\cap U$.
Then there exists a continuous family $G(t)\in\widetilde{\GL}_2(\R)$, $t\in\R$, such that $\sigma\cdot G(t)\to\overline{\sigma}$, for $t\to+\infty$, where $\overline{\sigma}\in\Theta_{\EEE}$ is a stability condition having $\AA_{\EEE}$ as heart and $\overline{\phi}(i_*\EE_2)=\overline{\phi}(i_*\EE_1[1])=\overline{\phi}(i_*\EE_0[2])=1$.
Hence $\overline{\sigma}\in\overline{U}$.
By Theorem \ref{thm:boundary}, there exist two stability conditions $\sigma_1\in\Theta_{\EEE}\cap\WW_{\EE_2}^+$ and $\sigma_2\in\WW_{\EE_0}^-$.
If $\phi_k$ denotes the phase function in $\sigma_k$ ($k=1,2$), we must have
\[
\phi_1(i_*\EE_0[2])>\phi_1(i_*\EE_2)>\phi_1(i_*\EE_1[1])
\]
and
\[
\phi_2(i_*\EE_1[1])>\phi_2(i_*\EE_0[2])>\phi_2(i_*\EE_2).
\]
But then $\sigma_1\in\Theta_{\EEE}\cap\Theta_{\tau_2^{\pm 1}\EEE}\cap\overline{U}$ and $\sigma_2\in\Theta_{\EEE}\cap\Theta_{\tau_1^{\pm 1}\EEE}\cap\overline{U}$.
Since the subsets $\Theta$ are open, this is enough to conclude that $\Theta_\EEE\cap\Theta_{\tau_j^{\pm 1}\EEE}\cap U\neq\emptyset$, for $j=1,2$.
This shows the claim.

In general, let $\sigma\in\Theta_\SSS\cap U$, for an ordered quivery collection $\SSS=\{S_0,S_1,S_2\}$.
Then, proceeding as above, $\Theta_\SSS\cap U$ is connected and we can find a stability condition $\overline{\sigma}$ in the closure of the $\widetilde{\GL}_2(\R)$-orbit of $\sigma$ such that $S_0$, $S_1$, and $S_2$ are $\overline{\sigma}$-stable of the same phase, that is $\overline{\sigma}\in\Theta_\SSS\cap\overline{U}$.
But then, by \eqref{eq:Utah29909}, $\overline{\sigma}\in\overline{\Theta}_{\EEE}$, for some exceptional collection $\EEE$ of vector bundles on $\P^2$. This gives $\Theta_\EEE\cap\Theta_\SSS\cap U\neq\emptyset$ and this intersection is connected, which completes the proof.
\end{Prf}

\medskip

We conclude the section by making a comparison between $\Stab_a$ and $\Stab^\dag(\DD_0)$.
To this end, we define $\Stab_g$ as the set of geometric stability conditions which, up to the action of $\C$, are of the form $\sigma_{a,b}$ with
\begin{align*}
&\Im a>0\\
& \Re b>-B\cdot\Re a+\frac{1}{2}B^2,
\end{align*}
where as in Definition \ref{def:setG}, $B:=-\frac{\Im b}{\Im a}$.
By Theorem \ref{thm:geom-stability}, all pairs $(a,b)\in\C^2$ satisfying the above inequalities are actual stability conditions.
This implies that $\Stab_g$ is an open, connected, and simply-connected subset of $\Stab(\DD_0)$.
Moreover, up to the action of $\widetilde{\GL}_2(\R)$, we can assume the central charge of a stability condition in $\Stab_g$ to take the form (see
\cite{Aaron-Daniele})
\[
Z^{t,m}(-)=-\int_{\P^2}e^{-(t+im)h}\ch(-),
\]
for $t,m\in\R$, $m>0$, and $h$ the class of a line in $\P^2$. In such a case, we denote the corresponding stability condition by $\sigma^{t,m}$.

Let $\EEE=\{\EE_0,\EE_1,\EE_2\}$ be an exceptional collection of vector bundles on $\P^2$. Set $\ch(\EE_j)=(r_j,d_j,c_j)$, $\mu_j:=d_j/r_j$, and $\Delta_j:=\frac{1}{2}\left(1-\frac{1}{r_j^2}\right)$ ($j=1,2,3$). Note that
$\mu_0 < \mu_1 < \mu_2 < \mu_0 + 3$.

\begin{Prop}\label{prop:geomvsalg}
We have $\sigma=\sigma^{t,m}\in\Theta_{\EEE}\cap\Stab_g$ only if $(t,m)$ is contained in the open semicircle with center $(C,0)$, where
\begin{equation}\label{eq:center}
C:=\frac{1}{2}(\mu_0+\mu_2)+\frac{\Delta_0-\Delta_2}{\mu_2-\mu_0}
\end{equation}
and radius $R:=\sqrt{\rho}$, where
\begin{equation}\label{eq:radius}
\rho:=\left(\frac{\Delta_0-\Delta_2}{\mu_2-\mu_0}\right)^2+\frac{1}{4}(\mu_2-\mu_0)^2-\left(\Delta_0+\Delta_2\right)>0.
\end{equation}
\end{Prop}
\begin{Prf}
First of all, let $\sigma=\sigma^{t,m}\in\Theta_{\EEE}\cap\Stab_g$.
Then, by Lemma \ref{lem:thetaP2},
\begin{enumerate}
\item $\phi(i_*\EE_0)<\phi(i_*\EE_1)<\phi(i_*\EE_2)$ and
\item \label{enum:tmineq2} $\phi(i_*\EE_0)+1=\phi(i_*\EE_0[1])<\phi(i_*\EE_2)$.
\end{enumerate}
As a consequence of (\ref{enum:tmineq2}), $(t,m)$ lies in the region bounded by
\[
\frac{\Im Z^{t,m}(i_*\EE_0)}{\Re Z^{t,m}(i_*\EE_0)}=\frac{\Im Z^{t,m}(i_*\EE_2)}{\Re Z^{t,m}(i_*\EE_2)},
\]
Making it explicit, we have
\[
m^2+\left(t-\frac{r_0c_2-r_2c_0}{d_2r_0-d_0r_2}\right)^2=-2\frac{d_0c_2-d_2c_0}{d_2r_0-d_0r_2}+\left(\frac{r_0c_2-r_2c_0}{d_2r_0-d_0r_2}\right)^2.
\]
As observed in Appendix \ref{app:DP},
\[
\frac{c_j}{r_j}=\frac{1}{2r^2}-\frac{1}{2}+\frac{\mu_j^2}{2}=-\Delta_j+\frac{\mu_j^2}{2},
\]
for $j=0,1,2$. Substituting we immediately deduce \eqref{eq:center} and \eqref{eq:radius}.
The fact that $\rho>0$ is again a straightforward computation, using
\[
0=\chi(\EE_2,\EE_0)=r_0r_2\left(1-\frac{3}{2}(\mu_2-\mu_0)+\frac{1}{2}(\mu_2-\mu_0)^2-(\Delta_0+\Delta_2)\right).
\]
\end{Prf}

Using Lemma \ref{lem:StabilitySkyscrapers} and a deformation argument it can be proved that the statement of the previous proposition is actually an \emph{if and only if}.

\begin{Rem}\label{rmk:Bonn1107}
By \cite[Proposition 5.1]{Goro-Ruda:Exceptional}, we have
\[
\frac{\Delta_0-\Delta_2}{\mu_2-\mu_0}=\frac{3}{2}\cdot\frac{r_0^2-r_2^2}{r_0^2+r_2^2+(c_2 r_0-c_0 r_2)^2}\in\left[-\frac{3}{2},\frac{3}{2}\right].
\]
Hence, if $m>3/\sqrt{2}$, then $\sigma^{t,m}\notin\Theta_{\EEE}$.
\end{Rem}

\section{Simply-connectedness}\label{sec:AlgebraicStability2}

We can now prove the simply-connectedness of $\Stab^\dag(\DD_0)$:

\begin{Thm}\label{thm:SimplyConn}
The connected component $\Stab^\dag(\DD_0)$ is simply-connected.
\end{Thm}

The idea of the proof is very simple: by an elementary topological argument, using what we proved in the previous section, we first reduce Theorem \ref{thm:SimplyConn} to proving that $\Stab_a$ is simply-connected.
To show this last assertion, we associate to every loop in $\Stab_a$ a word in the generators of the affine braid group
$B_3$. Then the simply-connectedness of $\Stab_a$ will be equivalent to the fact that $B_3$ acts freely on the set of
ordered quivery subcategories, and that, for every relation in $B_3$, we can find a corresponding loop that is contractible.

The main reason we involve Bridgeland's description of the set $\Stab_a$ is the following: the loci of degenerate
stability conditions appearing in Remark \ref{rem:degenerate} are rather implicit in our description of
$\Stab^\dag(\DD_0)$; however, they are essential for the simply-connectedness of the space.

\begin{Rem}\label{Rem:sc-slice}
Following \cite{Bridgeland:stab-CY}, let $\Stab_n \subset \Stab^\dag(\DD_0)$ be the subset of
normalized stability conditions with $Z(k(x)) = -1$. Denote by $\overline{\Stab}$ the quotient $\Stab^\dag(\DD_0)/\C$, which must also be simply-connected. By the results of the previous
section, there always exist semistable objects of class $[k(x)]$ in $\Stab^\dag(\DD_0)$; hence
$Z(k(x))$ is never zero. It follows that the subset $\Stab_n$ already surjects onto
$\overline{\Stab}$. This surjection is a Galois covering $\overline{\Stab} \cong \Stab_n/\Z$, where the
action by $n \in \Z$ is given as the shift $[2n]$; by the simply-connectedness of $\overline{\Stab}$, it
follows that $\Stab_n \cong \overline{\Stab} \times \Z$.

In particular, there is a connected component of ``very normalized'' stability conditions
$\Stab_{vn} \subset \Stab_n$ containing the geometric stability conditions where the skyscraper
sheaves are semistable of phase 1. It is a global slicing with respect to the $\C$-action, and
simply-connected. It is invariant under spherical twists and tensoring with line bundles (i.e.,
invariant under the subgroup $\Gamma_1(3) \subset \Aut \DD_0$ of Theorem \ref{thm:autoequiv-group}).
Presumably, Bridgeland's Conjecture 1.2 in \cite{Bridgeland:stab-CY} could be modified to use this
connected component $\Stab_{vn}$ rather than its open subset $\Stab_{n}^0(X)$ in the notation of
\cite{Bridgeland:stab-CY}.
\end{Rem}

\begin{Lem}\label{lem:VanKampen}
Let $X$ be a topological space such that
\[
X=A\cup\bigcup_{n\in I}B_n
\]
where $I$ is an arbitrary set of indices and
\begin{itemize}
\item $A$ and all $B_n$ are open, connected, and simply-connected;
\item $A\cap B_n$ is non-empty and connected, for all $n\in I$;
\item $B_n\cap B_m=\emptyset$, for $n\neq m$.
\end{itemize}
Then $X$ is simply connected.
\end{Lem}
\begin{Prf}
An inductive application of the classical Seifert--Van Kampen Theorem shows
that for all finite subsets $N\subset I$
\[
X_N:=A\cup\bigcup_{n\in N}B_n
\]
is connected and simply-connected. However, by compactness, any loop in $X$ is contained
in $X_N$ for some finite subset $N\subset I$. Hence it is contractible, as required.
\end{Prf}

To prove Theorem \ref{thm:SimplyConn}, we use the previous lemma with $A:=\Stab_a$, and the family $B_n$ as $\Phi(U)$, for $\Phi$ an autoequivalence of $\DD_0$ which belongs to the subgroup generated by spherical twists ${\ST}_\FF$ ($\FF$ an exceptional vector bundle on $\P^2$).
By Theorem \ref{thm:geom-stability}, $U \cong \C \times G$ is open, connected and simply-connected;
hence the same holds for all $\Phi(U)$. Also, we have $U\cap\Phi(U)=\emptyset$ unless $U=\Phi(U)$.

Thus  Theorem \ref{thm:SimplyConn} follows from Corollary \ref{cor:BoundaryIsAlgebraic}, Corollary
\ref{cor:IntersectionIsConnected}, and the following proposition:

\begin{Prop}\label{prop:main}
$\Stab_a$ is simply-connected.
\end{Prop}

Before proving Proposition \ref{prop:main}, we need a few lemmata.

\begin{Lem}\label{lem:firststabilityP2}
Let $\SSS=\{S_0,S_1,S_2\}$ and $\RRR=\{R_0,R_1,R_2\}$ be two ordered quivery collections.
Assume that $\Theta_\SSS\cap\Theta_\RRR\neq\emptyset$. Then, either $\Theta_{\SSS}=\Theta_{\RRR}$, or there exists a stability condition $\overline{\sigma}=(\overline{Z},\overline{\PP})\in\partial\Theta_{\SSS}\cap\Theta_{\RRR}$ such that the image of $\overline{Z}$ is contained in a line.
\end{Lem}
\begin{Prf}
By hypothesis, either $\Theta_{\SSS}=\Theta_{\RRR}$, or there exists a stability condition $\sigma\in\partial\Theta_{\SSS}\cap\Theta_{\RRR}$.
Now, we proceed as in the proof of Corollary \ref{cor:IntersectionIsConnected}: for every stability condition in $\partial\Theta_{\SSS}\cap\Theta_{\RRR}$ there exists a sequence $G_k\in\widetilde{\GL}_2(\R)$ ($k\in\N$) such that $\sigma\cdot G_k\to\overline{\sigma}$, where $\overline{\sigma}=(\overline{Z},\overline{\PP})$ is a stability condition in $\Theta_{\RRR}$ such that the image of $\overline{Z}$ is contained in a line. But then $\overline{\sigma}\in\partial\Theta_{\SSS}\cap\Theta_{\RRR}$, as wanted.
\end{Prf}

\begin{Lem}\label{lem:secondstabilityP2}
Let $\SSS=\{S_0,S_1,S_2\}$ and $\RRR=\{R_0,R_1,R_2\}$ be two ordered quivery collections.
If $\overline{\sigma}=(\overline{Z},\overline{\PP})\in\partial\Theta_{\SSS}\cap\Theta_{\RRR}$ is such that the image of
$\overline{Z}$ is contained in a line, then there exists $\gamma = \gamma_s\cdot\ldots\cdot\gamma_1\in B_3$,
$\gamma_k\in\{\tau_0^{\pm 1},\tau_1^{\pm 1},\tau_2^{\pm 1}\}$ for all $k\in\{1,\ldots,s\}$, such that, up to reordering, $\RRR=\gamma\SSS$, and there exist real numbers $0=a_0<a_1<\ldots<a_s<a_{s+1}=1$ and a continuous path $\alpha \colon [0,1]\to\Stab_a$ such that $\alpha([a_k,a_{k+1}))\subset\Theta_{\gamma_k\ldots\gamma_1\SSS}$ and $\alpha(1)=\overline{\sigma}$.
\end{Lem}
\begin{Prf}
First of all, if the image of $\overline{Z}$ is contained in a line, then the quivery collection is uniquely determined, up to reordering.
Then, given $\overline{\sigma}$, we can deform it slightly as in the proof of Proposition \ref{prop:BoundaryTheta}, case (a).
In this way we can find a non-trivial $\gamma_1\in\{\tau_0^{\pm 1},\tau_1^{\pm 1},\tau_2^{\pm 1}\}$ such that $\overline{\sigma}\in\overline{\Theta}_{\gamma_1\SSS}\cap\partial\Theta_{\SSS}$.
If $\overline{\sigma}\in\Theta_{\gamma_1\SSS}$, then $\RRR=\gamma_1\SSS$, up to the action of $r$, and the lemma is proved.
Otherwise, we can iterate the previous argument, by replacing $\SSS$ with $\gamma_1\SSS$.
This process terminates as in the proof of Proposition \ref{prop:BoundaryTheta}.
\end{Prf}

\begin{Lem}\label{lem:SimplConn1}
Let $\SSS$ be an ordered quivery collection.
Then, for all $\gamma\in\{\tau_0^{\pm 1},\tau_1^{\pm 1},\tau_2^{\pm 1}\}$, $\Theta_\SSS\cup\Theta_{\gamma\SSS}$ is simply-connected.
\end{Lem}
\begin{Prf}
For simplicity, we assume $\gamma=\tau_1$.
By Lemma \ref{lem:thetaP2} and the Seifert-Van Kampen Theorem, we only need to show that $\Theta_{\SSS}\cap\Theta_{\tau_1\SSS}$ is connected.
But, using Lemma \ref{lem:thetaP2} again, as well as the remark following it, we have
\[
\Theta_{\SSS}\cap\Theta_{\tau_1\SSS}=\CC_\SSS\cap\left\{(m_0,m_1,m_2,\phi_0,\phi_1,\phi_2)\in\R^6\colon \phi_1>\phi_0,\,\phi_1>\phi_2\right\},
\]
which is clearly connected.
\end{Prf}

\begin{Prf} (Proposition \ref{prop:main})
Take a continuous loop $\alpha \colon [0,1]\to\Stab_a$.
By using the previous lemmata, there exist real numbers $0=a_0<a_1<\ldots< a_m=1$, $m\in\N$, and ordered spherical collections
\[
\{\MMM_k=\{M_0^k,M_1^k,M_2^k\}\}_{k\in\{1,\ldots,m\}}
\]
with $\MMM_{k+1}$ obtained from $\MMM_k$ by an element $\gamma_{k+1}\in\{\tau_0^{\pm 1},\tau_1^{\pm 1},\tau_2^{\pm
1}\}$, such that, up to replacing $\alpha$ with an homotopic path, $\alpha([a_{k-1},a_k))\subset\Theta_k:=\Theta_{\MMM_k}$ for $k\in\{1,\ldots,m\}$ and
$\alpha(0)\in\Theta_m\cap\Theta_1$.
Thus we can assign a word $W(\alpha) = \gamma_m \ldots \gamma_1$ in the generators
of $B_3$ to every loop $\alpha$. Using Lemma \ref{lem:SimplConn1} we deduce that the homotopy class $[\alpha]$ of
$\alpha$ is determined by $W(\alpha)$, and that $[\alpha]$ is in fact determined by the element in the free group $L$
generated by $r, \tau_0, \tau_1, \tau_2$ associated to $W(\alpha)$.

Now assume more specifically that the stability conditions $(Z, \PP) = \alpha(0) = \alpha(1)$ is given by $\PP((0,1]) =
\AA_{\EEE_1}$, $Z(i_*\OO_{\P^2}) = Z(i_*\Omega_{\P^2}(1)[1]) = Z(i_*\OO_{\P^2}(-1)[2]) = - \frac 13$.  Then any heart
$\PP((\phi, \phi+1])$ is a shift of $\AA_{\EEE_1}$. As an ordered quivery collection is determined, up to reordering, by
its heart, we have
\[
r^j W(\alpha)(\MMM_1)=\MMM_1,
\]
for some $j$. Since the braid group $B_3$ acts freely on the set of ordered quivery subcategories (by \cite[Theorem 5.6]{Bridgeland:TStruct}), we have $r^jW(\alpha) = \id_{B_3}$ in $B_3$. Due to the description of $B_3$ in terms of generators and relations, it follows that we have an identity in $L$ of the form 
\[
W(\alpha) = r^{-j}(h_1R_1^{\pm 1}h_1^{-1})\cdots(h_sR_s^{\pm 1}h_s^{-1}),
\]
with $R_1, \dots, R_s \in\{r\tau_ir^{-1}\tau_{i+1}^{-1},\tau_i\tau_{i+1}\tau_i\tau_{i+1}^{-1}\tau_i^{-1}\tau_{i+1}^{-1},r^3\}$ and $h_1,\ldots,h_s\in L$ arbitrary elements.

By Lemma \ref{lem:SimplConn2}, loops with associated words $\tau_i\tau_{i+1}\tau_i\tau_{i+1}^{-1}\tau_i^{-1}\tau_{i+1}^{-1}$
(or its inverse) can be contracted in $\Stab_a$.  This implies that $\alpha$ can be contracted in general, and so
$\Stab_a$ is simply-connected.
\end{Prf}

\begin{Lem}\label{lem:SimplConn2}
Let $\alpha$ be a loop with word $W(\alpha) = \tau_i\tau_{i+1}\tau_i\tau_{i+1}^{-1}\tau_i^{-1}\tau_{i+1}^{-1}$.
Then $\alpha$ is contractible.
\end{Lem}
\begin{Prf}
We may assume $i=1$. We will say that a loop $\alpha$ ``runs through the regions $U_1, \ldots, U_m$'' for open subsets
$U_i \subset \Stab_a$ if there are
$0=a_0 < a_1 < \ldots <a_{m} < 1$ with $\alpha([a_{k-1}, a_k)) \subset U_k$ and $\alpha([a_m,1]) \subset U_1$.

By assumption, the loop $\alpha$ runs through the regions $\Theta_1, \dots, \Theta_6$ given
by $\Theta_k = \Theta_{\MMM_k}$ and
\begin{align*}
\MMM_1 &= \{S_0,S_1,S_2\} &
\MMM_4 &= \{\Gamma,{\ST}_{S_0}^{-1}S_1[1],S_0[2]\} \\
\MMM_2 &= \{S_0,{\ST}_{S_1}^{-1}S_2,S_1[1]\} &
\MMM_5 &= \{{\ST}_{S_0}^{-1}S_1,{\ST}_{S_0}^{-1}S_2,S_0[2]\} \\
\MMM_3 &= \{\Gamma,S_0[1],S_1[1]\} &
\MMM_6 &= \{{\ST}_{S_0}^{-1}S_1,S_0[1],S_2\}
\end{align*}
for some ordered quivery collection $S_0, S_1, S_2$, where $\Gamma={\ST}_{S_0}^{-1}{\ST}_{S_1}^{-1}S_2$.

First of all observe that, by Lemma \ref{lem:thetaP2},
\[
\Theta_1\cap\Theta_2\cap\Theta_3\neq\emptyset
\]
and it is homeomorphic to the locus in $\Theta_1$ given by those stability conditions having phases such that $\phi(S_1)<\phi(S_0)<\phi(S_2)$ and $\phi(S_0)<\phi({\ST}_{S_1}^{-1}S_2)$.
This implies that we can replace $\alpha$ by a loop, which we will denote again $\alpha$, such that
$\alpha$ runs though the regions $\Theta_1,\Theta_3, \Theta_4, \Theta_5,\Theta_6$. Repeating the same argument on
$\Theta_4\cap\Theta_5\cap\Theta_6$, we can replace it by a loop that runs through $\Theta_1, \Theta_3, \Theta_4,
\Theta_6$.

Let $t_1\in(0,1)$ be such that $\alpha([0,t_1))\subset\Theta_1$ and $\alpha(t_1)=(Z_1,\PP_1)\in\Theta_3\cap\partial\Theta_1$.
By Lemma \ref{lem:firststabilityP2}, we can assume the image of $Z_1$ to be contained in the real line.
In such a case, we have $\phi_1(S_0)=\phi_1(S_1)=0$ and $\phi_1(S_2)=1$.
At the same time, by definition, we have
\begin{equation}\label{eqn:Tianjin1}
S_0^{\oplus m}\to{\ST}_{S_0}^{-1}S_1\to S_1
\end{equation}
and so ${\ST}_{S_0}^{-1}S_1$ is semistable as well of phase $0$ whose Jordan-H\"older filtration is given by \eqref{eqn:Tianjin1}.
By Lemma \ref{lem:semitostable}, $\alpha(t_1)\in\Theta_3\cap\partial\Theta_6\cap\partial\Theta_1$.
In particular, $\Theta_3\cap\Theta_6\neq\emptyset$.

Let $t_4\in(t_1,1)$ be such that $\alpha(t_4)\in\Theta_1\cap\partial\Theta_6$ and $\alpha((t_4,1])\subset\Theta_1$.
By Lemma \ref{lem:thetaP2}, the intersection $\Theta_1\cap\partial\Theta_6$ is given as the
region
\[
\CC_{\MMM_1}\cap\left\{(m_0,m_1,m_2,\phi_0,\phi_1,\phi_2)\in\R^6\colon \phi_1=\phi_0,\,\phi_0\leq\phi_2\leq\phi_0+1\right\}.
\]
in particular, it is connected and simply-connected.
Since $\alpha(t_1)\in\overline{\Theta}_1\cap\partial\Theta_6$, we can replace $\alpha$ by a
homotopic path for which $\alpha(t_1)=\alpha(t_4)\in\Theta_3\cap\partial\Theta_6\cap\partial\Theta_1$.

What we proved so far is that our original loop $\alpha$ is homotopic to a loop which can be
decomposed as a loop $\alpha' = \alpha([0, t_1] \cup [t_4, 1])$ contained in
$\Theta_1\cup\Theta_3$ and another loop $\beta = \alpha([t_1,t_4])$ which runs through $\Theta_3,\Theta_4,\Theta_6$.

Now consider just the loop $\beta$.
Let $t_2\in(t_1,t_4)$ be such that $\beta(t_2)\in\Theta_6\cap\partial\Theta_4$ and
$t_3\in(t_2,t_4)$ such that $\beta(t_3)\in\Theta_4\cap\partial\Theta_3$. Arguing as above, we can
replace $\beta$ by a homotopic loop for which $\beta(t_2)=\beta(t_3)\in\Theta_6\cap\partial\Theta_3\cap\partial\Theta_4$.
Hence, the loop $\beta$ can be decomposed as a loop $\alpha'' =\beta([t_2, t_3])$ contained in
$\Theta_4\cup\Theta_6$ and another loop $\alpha''' = \beta([t_1, t_2] \cup [t_3, t_4])$ contained in $\Theta_3\cup\Theta_6$.

Summing up, to prove that $\alpha$ is contractible, we only need to prove that all regions $\Theta_1\cup\Theta_3$, $\Theta_4\cup\Theta_6$, and $\Theta_3\cup\Theta_6$ are simply-connected.
Again, by Lemma \ref{lem:thetaP2} and the Seifert--Van Kampen Theorem, it is sufficient to show that the intersections $\Theta_1\cap\Theta_3$, $\Theta_4\cap\Theta_6$, and $\Theta_3\cap\Theta_6$ are connected.
For $\Theta_1\cap\Theta_3$, observe that it corresponds to the locus in $\CC_{\MMM_1}$ in which $\Gamma$ is stable and $\phi(S_0),\phi(S_1)<\phi(\Gamma)$.
This can be proved to be connected by proceeding in a similar way as in the last part of the proof of Lemma \ref{lem:StabilitySkyscrapers}: in this situation we use the $\widetilde{\GL}_2(\R)$-action to fix the values of $Z(\Gamma)$ and $Z(S_1)$. The region $\Theta_1\cap\Theta_3$ is then, up to the action of $\widetilde{\GL}_2(\R)$, cut out by half-planes, and so it is connected.
The intersection $\Theta_4\cap\Theta_6$ is analogous.
Finally, $\Theta_3\cap\Theta_6$ corresponds to the locus in $\CC_{\MMM_3}$ in which $S_2$ is stable, $\phi(S_0[1])<\phi(S_1[1])$, and $\phi(S_2)<\phi({\ST}_{S_0}^{-1}S_1)$, which is again connected by a similar argument. This completes the proof of the lemma.
\end{Prf}

\section{Group of autoequivalences}\label{sec:autoequivalences}

Let $\Aut \DD_0$ be the group of autoequivalences of $\DD_0$ up to isomorphism of functors, and let
$\Aut^\dag \DD_0$ be the subgroup of $\Aut \DD_0$ preserving the connected component
$\Stab^\dag(\DD_0)$.

The numerical $K$-group of $\DD_0$ is
$K(\DD_0)/K^\perp = K(\DD_0)/ \Z\cdot[k(x)] \cong \Z^{\oplus 2}$. Since the Euler
form is skew-symmetric, there is a natural map 
\begin{equation} 			\label{eq:AuttoSL2Z}
\Aut \DD_0 \to \SL(2, \Z)		\end{equation}
given by sending
an autoequivalence to its induced action on the numerical $K$-group.
Crucial for us will be the congruence subgroup
$\Gamma_1(3) \subset \SL(2, \Z)$ of matrices 
\[ \begin{pmatrix} a & b \\ c & d \end{pmatrix} \equiv
 \begin{pmatrix} 1 & b \\ 0 & 1 \end{pmatrix} 
\pmod{3} \]
It has generators $T = \begin{pmatrix} 1 & 1 \\ 0 & 1 \end{pmatrix}$
and $S = \begin{pmatrix} 1 & 0 \\ -3 & 1 \end{pmatrix}$
with a single relation given by $(ST)^3 = 1$.

As in the introduction, we denote by $\hat X$ the formal completion of $X$ along $\P^2$.

\begin{Thm}\label{thm:autoequiv-group}
\[ \Aut^\dag \DD_0 \cong \Z \times \Gamma_1(3) \times \Aut(\hat X). \]
\end{Thm}

We start by identifying the subgroup $\Gamma_1(3) \subset \Aut^\dag(\DD_0)$:
As observed in \cite[Section 7.3.5]{Aspinwall:Dbranes-CY}, there is a relation
\begin{equation}\label{eq:Z3Z-relation}
 \bigl({\ST}_{\OO_{\P^2}} \circ (\blank \otimes \pi^*\OO(1))\bigr)^3 \cong \Id.
\end{equation}
Due to the description of $\Gamma_1(3)$ by generators and relations, this induces a map
$\Gamma_1(3) \to \Aut^\dag \DD_0$. If we choose $([\OO_{\P^2}], [\OO_l])$ as a basis
of $K(\DD_0)/K^\perp$, then the composition with \eqref{eq:AuttoSL2Z} maps the generators to $S$ and
$T$, respectively; hence the composition
$\Gamma_1(3) \to \Aut^\dag \DD_0 \to \SL_2(\Z)$ is the standard inclusion of $\Gamma_1(3)$ as the
congruence subgroup given above. In particular, the action on $\DD_0$ is faithful, and we obtain:

\begin{Prop}				\label{prop:gamma1}
The subgroup of $\Aut^\dag \DD_0$ generated by
${\ST}_{\OO_{\P^2}}$ and $\blank \otimes \pi^*\OO(1)$
is isomorphic to $\Gamma_1(3)$.
\end{Prop}

Alternatively, one can prove that the composition
${\ST}_{\OO_{\P^2}} \circ (\blank \otimes \pi^*\OO(1))$ is isomorphic to the generator of
the natural $\Z_3$-action on $\Db_0([\C^3/\Z_3])$ via the derived McKay correspondence of
\cite{Mukai-McKay} (where $\Db_0([\C^3/\Z_3])$ denotes the bounded derived category of
$\Z_3$-equivariant coherent sheaves on $\C^3$ supported at the origin).

\begin{Lem}\label{lem:ST-in-gamma13}
For any exceptional vector bundle $\EE$ on $\P^2$, its associated spherical twist ${\ST}_\EE$ is
contained in the subgroup $\Gamma_1(3)$ generated by ${\ST}_{\OO_{\P^2}}$ and
$\blank \otimes \pi^* \OO(1)$.
\end{Lem}
\begin{Prf}
Since, for all autoequivalences $\Phi\in\Aut(\DD_0)$ and for all spherical objects $F\in\DD_0$,
\[
\Phi\circ{\ST}_{F}\circ\Phi^{-1}\cong{\ST}_{\Phi(F)},
\]
it is sufficient to show that there exists $g \in \Gamma_1(3)$ with $i_*\EE \cong g\left(\OO_{\P^2}\right)$.
By \cite{Goro-Ruda:Exceptional}, $\EE$ is contained in a mutation of the exceptional collection $\EEE_1$.
Remark \ref{rmk:GoroRuda} completes the proof.
\end{Prf}

Restricting to the subgroup $\Aut^\dag \DD_0$ makes it possible to control autoequivalences via the
following proposition:

\begin{Prop} \label{prop:geom-autos}
Let $\Phi \in \Aut \DD_0$ be an autoequivalence such that there exist two geometric stability
conditions $\sigma, \sigma'$ with $\Phi(\sigma) = \sigma'$. Then $\Phi$ is isomorphic to the
composition of an automorphism of $\hat X$ with $\blank \otimes \OO(n)[k]$.
\end{Prop}
\begin{Prf}
From the description of geometric stability conditions it follows
that $\Phi$ sends skyscraper sheaves $k(x)$, $x \in \P^2$, to shifts of skyscraper sheaves. More precisely, after replacing 
$\Phi$ by $\Phi \circ \blank[k]$ for some $k \in \Z$, we may assume that
for every $x \in \P^2$ there is $x' \in \P^2$ with
$\Phi(k(x)) \cong k(x')$.

The proposition follows now directly from the recent results in \cite{Lunts-Orlov:uniqueness,
alberto-paolo:FM-supported}.
Since the argument is quite standard, we give only a brief sketch.
By \cite[Theorem 1.1]{alberto-paolo:FM-supported}, every autoequivalence of $\DD_0$ is of \emph{Fourier-Mukai type}, i.e., there exists an object $\UU\in\Db(\Qcoh_{\P^2\times\P^2}(X\times X))$ such that $\Phi\cong\Phi_\UU:=(p_1)_*\left(\UU\otimes p_2^*(\blank)\right)$, where $\Qcoh_{\P^2\times\P^2}(X\times X)$ is the category of quasi-coherent sheaves on $X\times X$ supported on $\P^2\times\P^2$, all functors are supposed to be derived, and $p_1,p_2$ denote the two projections.

Let $X_n$ be the $n$-th infinitesimal neighborhood of $\P^2$ inside $X$.
By using \cite[Lemma 4.3]{Bridgeland:EqFMT}, we can show that $\UU$ is actually a sheaf and it maps $\Coh_0$ (resp.\ $\Coh X_n$) to $\Coh_0$ (resp.\ $\Coh X_n$).
By tensoring with a line bundle on $X$, we can assume that $\Phi(\OO_{X_n})\cong\OO_{X_n}$, for all $n$.
Hence, arguing as in \cite[Corollary 5.23]{Huybrechts:FM}), since $\Phi_{\UU}(k(x))\cong k(x')$, there exists a family of compatible automorphisms $u_n\colon X_n\to X_n$, $n\in\Z_{\geq 0}$, which induces an automorphism $u\colon \widehat{X}\to\widehat{X}$ and such that $\Phi_{\UU}\cong u^*$, as wanted.
\end{Prf}

\begin{Prf} (Theorem \ref{thm:autoequiv-group})
Given $\Phi\in\Aut^\dag \DD_0$, pick an arbitrary geometric stability condition $\sigma \in U$.
By Corollary \ref{cor:ConnectedComponent}, there exists a stability condition
$\sigma' \in \overline{U}$ and a composition $\Psi$ of spherical twists associated to exceptional
vector bundles with $\Psi \circ \Phi(\sigma) = \sigma'$. The stability condition
$\Psi \circ \Phi(\sigma)$ has no semistable objects of class $[k(x)]$; thus actually $\sigma' \in
U$. By Proposition \ref{prop:geom-autos} and Lemma \ref{lem:ST-in-gamma13}, $\Phi$ is contained in
the group generated by $\Gamma_1(3)$, shifts, and $\Aut(\hat X)$.

As the actions by $\Z$, by $\Gamma_1(3)$, and by $\Aut(\hat X)$ commute, we get a surjective map
\[ \Z \times \Gamma_1(3) \times \Aut(\hat X) \to \Aut^\dag \DD_0.\]
The restriction 
$\Aut(\hat X) \to \Aut(\DD_0)$ is injective by \cite[Theorem 1.1]{alberto-paolo:FM-supported}.
Since $\Z \times \Aut(\hat X)$ acts by $\pm 1$ on $K/K^\perp$, its intersection with $\Gamma_1(3)$
is trivial, and the above map is an isomorphism.
\end{Prf}

\section{$\Pi$-stability and Global Mirror Symmetry}\label{sec:MS}

In this section, we outline how our results fit into expectations coming from mirror symmetry for
the local $\P^2$. Mirror symmetry for the local $\P^2$ has been discussed in many places of the
mathematical physics literature, see e.g. \cite{AGM:measuring, Diaconescu-Goms:fractional}; our presentation
follows \cite{ABK:topological_strings} and \cite{Aspinwall:Dbranes-CY} most closely.

\subsection{Monodromy and autoequivalences}
The family of mirror partners to the local $\P^2$ can be constructed
explicitly from the following family of genus one curves:
The equation
\[ X_0^3 + X_1^3 + X_2^3 - 3 \psi X_0 X_1 X_2 = 0 \]
cuts out a surface $S \subset \P^2 \times \C$. At $\psi^3 = 1$ and $\psi = \infty$, the fibers
are singular; all other fibers of $S$ over $\C$ are smooth genus one curves.
There is a
$\mu_3$-action on $S$ given by $X_0 \mapsto \omega^{-1} X_0$ and $\psi \mapsto \omega \psi$, and
leaving the other variables invariant, where $\omega=\exp(2\pi i/3)$.
Let $\YY$ be the quotient 
\[ \bigl(S \setminus \{ \psi^3=1 \}\bigr)/\mu_3 \]
of the union of the smooth fibers by the group action; then $\YY$ is a family of smooth elliptic
curves over $(\C - \mu_3)/\mu_3$.

In fact, the base is the moduli space $\MG13 \cong (\C - \mu_3)/\mu_3$ of elliptic curves with
$\Gamma_1(3)$-level structure.  We can also think of $\MG13$ as $\P^1$ with the points
$z = -\frac 1{27}$ and $z= 0$ removed, and a
stacky $\Z_3$ point at $z= \infty$ (where we set $z = - \frac 1{(3\psi)^3}$).  The fundamental group
of $\MG13$ is $\Gamma_1(3)$. It is generated by the loops $\gamma_{-\frac 1 {27}}, \gamma_0$ around $-\frac 1{27}$ and
$0$; as their composition is a loop around $z = \infty$, they satisfy
$\bigl(\gamma_{-\frac 1{27}} \gamma_0 \bigr)^3 = 1$.

Given any $z \in \MG13$, one can determine a basis of first homology $H_1(\YY_z)$ of the fibers by
choosing a path from $z$ to $-\frac 1{27}$ and from $z$ to $\infty$; the basis is then given by the two
corresponding vanishing cycles $\bar A_z$ and $\bar B_z$. This basis yields an identification of $\pi_1(\MG13)$ as a
subgroup of $SL_2(\Z)$ by its monodromy action on the first homology $H_1$ of the fibers of $\YY$.
Explicitly, we get
\[ \gamma_{-\frac 1{27}} \equiv \begin{pmatrix} 1 & 0 \\ -3 & 1 \end{pmatrix}, \quad
\gamma_0 \equiv \begin{pmatrix} 1 & 1 \\ 0 & 1 \end{pmatrix}.
\]

A well-known principle of mirror symmetry states that monodromies in the mirror family $\YY$ lift to
autoequivalences in the derived category $\DD_0$: it is implied by homological mirror symmetry and
has been applied and verified e.g. in \cite{Seidel-Thomas:braid, Horja:autoequivalences}. Theorem
\ref{thm:autoequiv-group} gives another incarnation of this principle, as the action of $\Gamma_1(3)
\subset \Aut \DD_0$ on $K(\DD_0)/K^{\perp} \cong \Z^{\oplus 2}$ matches the action of
$\Gamma_1(3) \cong \pi_1(\MG13)$
if we identify ${\ST}_\OO$ with $\gamma_{-\frac 1{27}}$ and $\blank \otimes \OO(1)$ with $\gamma_0$.

\subsection{Period integrals and $\Stab^\dag(\DD_0)$}
However, in the spirit of \cite{Bridgeland:spaces}, there is also more geometric connection between the
mirror moduli space $\MG13$ and the space of stability conditions.
The periods in this mirror construction are given by integrals over the meromorphic differential
form $\lambda = \ln \frac{X_2}{X_3} \frac{dX_1}{X_1}$. More precisely, let $\YY^0 \subset \YY$ be
the complement of the set of poles of $\lambda$, and $\tilde \YY^0$ be the cover on which $\ln
\frac{X_2}{X_3} $ is well-defined. Following \cite{ABK:topological_strings}, one can choose 
a family of cycles $A_z, B_z \in H_1(\tilde \YY^0_z)$ that project to $\bar A_z, \bar B_z \in H_1(\YY_z)$, and a third
family of cycles $C_z \in H_1(\tilde \YY^0_z)$ and define the period integrals as:

\[ \Pi(z) = \begin{pmatrix} \int_B \lambda \\ \int_A \lambda \\ \int_C \lambda \end{pmatrix} \]

The authors show that if $A, B, C$ are chosen appropriately, then the action of
$\pi_1(\MG13) \cong \Gamma_1(3)$ on these 3 periods matches the action of $\Gamma_1(3) \subset
\Aut \DD_0$ on $K(\DD_0) \cong \Z^{\oplus 3}$.

We will now ignore the construction of period integrals and instead just consider their Picard-Fuchs
equation; with $\theta_z := z \frac{d}{dz}$ it is given by
\begin{equation} \label{eq:Picard-Fuchs}
 \bigl( \theta_z^3  + 3z  \theta_z (3\theta_z + 1)(3\theta_z + 2) \bigr) \Pi = 0,
\end{equation}
and has singularities at $z = 0$, $z = - \frac 1{27}$ and $z = \infty$.

Using an \emph{Ansatz} and solving for the coefficients of the power series,
one can find expansions of three linearly independent solutions around
$z = 0$ and $\psi = 0$, respectively (see also \cite[Section 6]{ABK:topological_strings} and
\cite[Section 7.3]{Aspinwall:Dbranes-CY}).

Around $z = 0$, we make the standard branch choice of $\ln z$ for $z \in \C \setminus \R_{\le 0}$,
and get as expansions (compare with \cite[Section 6.2]{ABK:topological_strings})
\begin{align*}
\omega_0(z) &= 1 \\
\omega_1(z) &= \frac 1{2\pi i} \left(\ln z 
		+3 \sum_{n=1}^\infty \frac{(3n-1)!}{n!^3} (-z)^n \right) \\
\omega_2(z) &= \frac 1{(2 \pi i)^2} \left(
			(\ln z)^2 
			+ 6 \ln z \cdot \sum_{n=1}^\infty \frac{(3n-1)!}{n!^3} (-z)^n
			+ \sum_{n=1}^\infty l_n z^n \right)
\end{align*}
where the differential equation defines the $l_n$ recursively:
\begin{multline*}
 l_n = -\frac 1{n^3}  \left( (3n-1)(3n-2)(3n-3) l_{n-1} 
				+ 18 \cdot (-1)^n \frac{(3n-1)!}{n!^3} \cdot n^2 \right. \\ \left.
				- 18 \cdot (-1)^n \frac{(3n-4)!}{(n-1)!^3} \left(27 n^2 - 36 n + 11\right)
			\right)
\end{multline*}
Similarly, the power series expansion of a basis of solutions nearby $z = \infty$ are given by 
$\varpi_0(\psi) = 1$ and:
\begin{align*}
\varpi_1(\psi) &= \frac 1{2\pi i} \sum_{n=1 \atop 3\nmid n}^\infty
	\frac{\Gamma\left(\frac n3 \right)}
		{\Gamma(n+1) \Gamma\left(1 - \frac n3 \right)^2} (3 \psi)^n \\
\varpi_2(\psi) &= \frac 1{2\pi i} \sum_{n=1 \atop 3\nmid n}^\infty
	\frac{\Gamma\left(\frac n3 \right)}
		{\Gamma(n+1) \Gamma\left(1 - \frac n3 \right)^2} (3 e^{\frac{2\pi i}3} \psi)^n 
\end{align*}
Here we use $\psi = - \frac 1{3\sqrt[3]{z}}$ with the branch choice $\frac{2\pi}3 < \arg \psi < \frac{4\pi}3$ for
$\abs{\arg z} < \pi$. 

Following Aspinwall, we define the solutions $a(z), b(z)$ of \eqref{eq:Picard-Fuchs} for $z \in \C^* \setminus \R_{<0}$
by setting 
\begin{equation} \label{eq:ab-exp0}
 a(z) = \omega_1(z) - \frac 12, \quad b(z) =  - \frac 12 \omega_2(z) + \frac 12 \omega_1(z) - \frac 14
\end{equation}
for $\abs{z} < \frac 1{27}$ and analytic continuation.
This analytic continuation is computed explicitly in \cite[Eqn. (286)]{Aspinwall:Dbranes-CY} and
\cite[Eqn. (6.22)]{ABK:topological_strings},
and gives the following expansion of $a(z)$ and $b(z)$ around $z = \infty$:
\begin{equation} \label{eq:ab-expinf}
a(z) = \varpi_1(z)-\frac{1}{2}, \quad b(z) = \frac{1}{3}\left(\varpi_1(z)-\varpi_2(z)-1\right)
\end{equation}

\begin{Thm}\label{thm:PowerSeries}
Fix a universal cover $\tildeMG13  \to \MG13$ together with its $\Gamma_1(3)$-action of deck transformations
and choose a fundamental domain $D\subset\tildeMG13$ that projects isomorphically onto $\C^* \setminus \R_{<0}$.
Then there is an embedding $I \colon \tildeMG13 \to \Stab^\dag(\DD_0)$ defined by the following properties:
\begin{enumerate}
\item For $(Z(z), \PP(z)) = I(z)$, the central charge is given by
\[ Z(E)(z) = -c(E) + a(z) \cdot d(E) + b(z) \cdot r(E), \]
for all $E\in K(\DD_0)$ (where we identify $a(z), b(z)$ with their analytic continuations from $D$ to $\tildeMG13$).
\item For $z \in D$, the stability condition $I(z)$ is geometric with $k(x)$ having phase 1.
\end{enumerate}
\end{Thm}

On the boundary of $I(D)$ we have several interesting special points: the point $z=0$ corresponds to the \emph{large
volume limit point}, where $\Im(a)\to+\infty$ and the central charge is approximately given by $Z(E) = - \int_{\P^2}
\ch(E) e^{-ah}$, where $h$ is the class of a line in $\P^2$; the limit stability condition as $z \to 0$ can be described
as a polynomial stability condition of \cite{large-volume}.  The point $z=\infty$ is the \emph{orbifold point}: the
heart of the bounded $t$-structure is $\AA_1 \cong \Coh_0 [\C^3/\Z_3]$ and the three simple objects of $\AA_1$ have the
same central charge equal to $-1/3$; this point is fixed under the $\Z_3$-action on $\DD_0$ given by relation
\eqref{eq:Z3Z-relation} (i.e., by tensor product in $\Coh_0 [\C^3/\Z_3]$ with a non-trivial one-dimensional
$\Z_3$-representation).  Finally, when $\psi=\omega$ (resp.\ $\psi=\omega^2$) and so $z=-\frac 1{27}$ (these are called
\emph{conifold points}), we have a singularity: indeed, $Z(\OO_{\P^2})=0$, resp.\ $Z(\OO_{\P^2}(-1))=0$, depending on
whether we approach $-\frac 1{27}$ from above or below.

The proof of the theorem is based on the following two observations:

\begin{enumerate}
\item \label{obs:inequality}
For all $z \in \C \setminus \R_{\le 0}$, the complex numbers $a(z), b(z)$ satisfy the inequalities of
Definition \ref{def:setG}.
\item \label{obs:monodromy}
The monodromy action of $\Gamma_1(3)$ on the solutions $a(z),b(z)$
(computed, for example, in \cite{Aspinwall:Dbranes-CY}) is compatible with the $\Gamma_1(3)$-action on
$\Stab^\dag$.
\end{enumerate}

We first verified Observation (\ref{obs:inequality}) by explicit computations using the computer algebra
package SAGE \cite{sage}.\footnote{The program used to test the inequalities is available for download from the authors'
homepages. It implements the power series expansion around $z=0$ and $\psi=0$ and tests the inequalities for
random complex numbers in their respective convergence domains.}
A complete argument is sketched in Appendix \ref{app:ineq}.
To prove Theorem \ref{thm:PowerSeries} we only need to show Observation \eqref{obs:monodromy}:

\begin{Prf} (Theorem \ref{thm:PowerSeries})
By  Observation \eqref{obs:inequality} and Theorem \ref{thm:geom-stability}, we obtain an embedding
$I\colon D \into U \subset \Stab^\dag(\DD_0)$.
By Bridgeland's deformation result, the extension of $I$ to $\tildeMG13$ is unique, if it exists.

Now, we can extend $I$ to the $\Gamma_1(3)$-translates of $D$ uniquely by requiring it to be
$\Gamma_1(3)$-equivariant. Hence, it remains to check that this extension of $I$ glues along the translates
of $\partial D \subset \tildeMG13$, and is compatible with the requirement that $a(z), b(z)$ are solutions to the
Picard-Fuchs equation.

Let $\gamma_0$ be the loop going in positive direction around the origin $z = 0$, and
$\gamma_\infty$ the loop around $z= \infty$ acting on $\psi$ by
$\psi \mapsto e^{\frac{2\pi i}3} \psi$. Then, by the $\Gamma_1(3)$-equivariance, it is in fact enough to check
the glueing along $\overline{D} \cap \gamma_0(\overline{D})$ lying above $(-\frac 1{27}, 0) \subset \C$ 
in the $z$-plane, and along $\overline{D} \cap \gamma_\infty(\overline{D})$, lying above the line segments
$(0,1) \cdot e^{\frac 23 \pi i}$ and $(0,1) \cdot e^{\frac 43 \pi i}$ in the $\psi$-plane.

The action of $\gamma_0$ on the solutions is given by
\[
\omega_1(z)  \mapsto \omega_1(z) + 1, \quad \omega_2(z)  \mapsto \omega_2(z) + 2 \omega_1(z) + 1.
\]
The action of $\blank \otimes \OO(1)$ on the set of geometric stability conditions
$\sigma_{a, b}$ of Theorem \ref{thm:geom-stability} is given by $a \mapsto a + 1$ and
$b \mapsto b - a - \frac 12$. Using the expansions in equation \eqref{eq:ab-exp0}, we see that
the induced action of $\gamma_0$ on $a(z)$ and $b(z)$ matches exactly; hence the definition
of $a(z), b(z)$ on $\gamma_0(D)$ by analytic continuation agrees with the implicit definition given by the requirement
that $I$ is $\gamma_0$-equivariant; on the other hand, when $\arg(z) = \pi$ we have $B = 0$ and
$b(-\frac 1{27}) = 0$ (see Appendix \ref{app:ineq}), and it follows that 
$a(z), b(z)$ still satisfy the inequalities of Definition
\ref{def:setG} for $z \in \overline{D} \cap \gamma_0(\overline{D})$, i.e., for $z = (-\frac 1{27}, 0)$
with $\arg(z) = \pi$. Then Theorem \ref{thm:geom-stability} implies that $I$ glues along this boundary component of $D$ 
within the geometric chamber.

Similarly, the action of $\gamma_\infty$ on the space of solutions is computed in terms of the 
expansions around $z = \infty$ as
\[
\varpi_1(z) \mapsto \varpi_2(z), \quad 
\varpi_2(z) \mapsto - \varpi_1(z) - \varpi_2(z).
\]
The central charges of the three simple objects in the quiver category $\AA_1$ are given by
\begin{align*}
Z(\OO_{\P^2})(z) &= b(z) = \tfrac{1}{3}\varpi_1(z) - \tfrac 13 \varpi_2(z) - \tfrac 13 \\
Z(\Omega_{\P^2}(1)[1])(z) &= -2 b(z) + a(z) - \tfrac 12
= \tfrac{1}{3}\varpi_1(z) + \tfrac 23 \varpi_2(z) - \tfrac 13 \\
Z(\OO_{\P^2}(-1)[2])(z) &= b(z) - a(z) - \tfrac 12 = -\tfrac{2}{3}\varpi_1(z) - \tfrac 13 \varpi_2(z) - \tfrac 13
\end{align*}
The autoequivalence $\bigl({\ST}_{\OO_{\P^2}} \circ (\blank \otimes \pi^*\OO(1))\bigr)^{-1}$
permutes these 3 objects and preserves the heart of the t-structure $\AA_1$; hence it is easy
to see that its action on the central charge matches the monodromy $\gamma_\infty$. 
\end{Prf}

\appendix

\section{Bounds on stable Chern classes after Dr\'ezet-Le Potier}\label{app:DP}

We give a brief review and a reformulation of the main result of \cite{Drezet-LePotier}.
We recall that for a torsion-free sheaf $\FF$ on $\P^2$, its \emph{slope} is defined by
$\mu(\FF) = \frac{d(\FF)}{r(\FF)}$, giving the following notion of stability: 
\begin{Def} \label{def:stablesheaf}
A torsion-free sheaf $\FF$ on $\P^2$ is called \emph{slope-stable} if the inequality $\mu(\FF') < \mu(\FF)$ holds for all saturated subsheaves
$\FF' \subset \FF$.
\end{Def}

The \emph{discriminant} of $\FF$ is defined by $\Delta(\FF) = \frac{d(\FF)^2}{2r(\FF)^2} - \frac{c(\FF)}{r(\FF)}$.

Let $\AAA$ be the set of all $\alpha \in \Q$ such that there exists
an exceptional vector bundle on $\P^2$ with slope $\alpha$. For any 
$\alpha \in \AAA$, we call its rank $r_\alpha$ the smallest integer $r > 0$
such that $r\alpha \in \Z$. We call
$\Delta_\alpha := \frac 12 \left(1 - \frac 1{r_\alpha^2}\right)$
its discriminant.
It follows from Riemann-Roch that the rank and discriminant of an exceptional
vector bundle with slope $\alpha$ (if it exists) are uniquely determined by these formulas.
Similarly, any non-exceptional stable sheaf satisfies $\Delta \ge \frac 12$.

For two rational numbers with $3 + \alpha + \beta \neq 0$, Dr\'ezet and Le
Potier define the operation
\[ \alpha . \beta := \frac{\alpha + \beta}2
	+ \frac{\Delta_\beta - \Delta_\alpha}{3 + \alpha - \beta}
\]

Let $\DDD$ be the set of rational numbers of the form $\frac{p}{2^q}$ for
$p \in \Z, q \in \Z_{\ge 0}$. One defines a function
$\epsilon \colon \DDD \to \Q$ inductively by
$\epsilon(n) := n$ for $n \in \Z$ and
\[
\epsilon\left(\frac{2p+1}{2^{q+1}}\right) :=
	\epsilon\left(\frac p{2^q}\right) .
\epsilon\left(\frac{p+1}{2^q}\right)
\]

\begin{Thm}[{\cite[Th\'eor\`eme A and chapitre 5]{Drezet-LePotier}}]
\label{thm:DP-A}
The set $\AAA$ of exceptional slopes is equal to the image $\epsilon(\DDD)$, and for each
slope $\alpha \in \AAA$ the exceptional vector bundle of slope $\alpha$ is unique.
\end{Thm}

Now define 
\begin{align*}
P(X) & := 1 + \frac 32 X + \frac 12 X^2 \\
p(x) & := \begin{cases}
P(-\abs{x}) & \abs{x} < 3 \\
0 				    & \text{otherwise} \end{cases} \\
\intertext{and, for any $\alpha \in \AAA$,}
p_\alpha(x) &:= p(x - \alpha) - \Delta_\alpha.
\end{align*}

If $\alpha, \beta$ are of the form given in Theorem \ref{thm:DP-A},
then $p_\alpha$ and $p_\beta$ are monotone decreasing and increasing,
respectively; they intersect in the point 
$(\alpha . \beta, \Delta_{\alpha . \beta})$.

\begin{Thm}[{\cite{Drezet-LePotier}, \cite[Theorem 16.2.1]{LePotier}}] 		\label{thm:DP-C}
Given an integer $r > 0$ and rationals $\mu, \Delta \in \Q$, there exists a
slope-stable sheaf $\EE$ on $\P^2$ with rank $r$, 
slope $\mu$, and discriminant $\Delta$ if and only if
\begin{enumerate}
\item \label{divisibility}
$r\mu \in \Z$ and $r(P(\mu) - \Delta) \in \Z$, and
\item For every $\alpha \in \AAA$ with $r_\alpha < r$ and $\abs{\alpha - \mu}
<3$, we have $\Delta \ge p_\alpha(\mu)$.
\end{enumerate}
\end{Thm}
The only change compared to \cite[Theorem 16.2.1]{LePotier} is that we replaced
[Gieseker-]stable with slope-stable, which is justified by \cite[Th\'eor\`eme
(4.11)]{Drezet-LePotier}.

This leads us to define (cf. \cite[Sect.\ 16.4]{LePotier})
$\delta_\infty^{DP} \colon \R \to [1/2,1]$ as
\[
\delta_\infty^{DP}:= \sup \stv{p_\alpha}{\alpha \in \AAA}.
\]
The necessary and sufficient condition for the existence of non-exceptional slope-stable sheaves can then be written as
\[
\Delta\geq\delta_\infty^{DP}(\mu).
\]
As discussed in \cite[Sect.\ 16.4]{LePotier}, this is equivalent to the formulation in 
Thereom \ref{thm:DP-C} for purely arithmetic reasons.

The first part of Theorem \ref{thm:DP} now follows immediately.
For the last assertion, let $(\mu_n, \Delta_n)$ be a sequence of distinct points in $S_E$ that
converges in $\R^2$ to $(\mu, \Delta)$.
For every $\alpha \in \AAA$, we have $\Delta_n \ge p_\alpha(\mu_n)$
for all $n \gg 0$ (in fact, this holds whenever $\alpha \neq \mu_n$).
By continuity, $\Delta \ge p_\alpha(\mu)$, and thus
$\Delta \ge \delta_\infty^{DP}(\mu)$, i.e., the accumulation point
$(\mu, \Delta)$ is contained in $S_\infty$.

\section{Bridgeland's stability conditions}\label{app:BridgelandFramework}

In this section we give a brief review of stability conditions on derived categories, following \cite{Bridgeland:Stab}.

Let $\DD$ be a triangulated category with good properties, e.g.\ the bounded derived category of coherent sheaves on a smooth and projective variety or $\DD_0$.
A \emph{stability condition} $\sigma$ on $\DD$ consists in a pair $(Z,\PP)$, where $Z: K(\DD)\to\C$ (\emph{central charge}) is an additive map and $\PP(\phi)\subset\DD$ are full, additive subcategories ($\phi\in\R$) satisfying:
\begin{enumerate}
\item for any $0 \neq E\in\PP(\phi)$ we have $Z(E)\neq0$ and $Z(E)/|Z(E)|=\exp(i\pi\phi)$;
\item $\forall\phi\in\R$, $\PP(\phi+1)=\PP(\phi)[1]$;
\item if $\phi_1>\phi_2$ and $A_j\in\PP(\phi_j)$, $j=1,2$, then $\Hom(A_1,A_2)=0$;
\item \label{enum:HN-filt}
for any $E\in\DD$ there is a sequence of real numbers $\phi_1>\dots >\phi_n$ and a collection of triangles $E_{j-1}\to E_j\to A_j$ with $E_0=0$, $E_n=E$ and
$A_j\in\PP(\phi_j)$ for all $j$.
\end{enumerate}

The collection of exact triangles in (\ref{enum:HN-filt}) is called the \emph{Harder-Narasimhan filtration} of $E$.
Each subcategory $\PP(\phi)$ is extension-closed and abelian. Its nonzero objects are said to be \emph{semistable} of phase $\phi$ in $\sigma$, and the simple objects (i.e., objects without proper subobjects or quotients) are said to be \emph{stable}.

For any interval $I\subset\R$, $\PP(I)$ is defined to be the extension-closed subcategory of $\DD$ generated by the subcategories $\PP(\phi)$, for $\phi\in I$. Bridgeland proved that, for all $\phi\in\R$, $\PP((\phi,\phi+1])$ is the heart of a bounded $t$-structure on $\DD$.
The category $\PP((0, 1])$ is called the \emph{heart} of $\sigma$.

\begin{Rem}\label{rmk:tstruct}
Let $\H:=\{z\in\C\colon z=|z|\exp(i\pi\phi),\,0<\phi\leq1\}$. If $\AA\subset\DD$ is the heart of a bounded $t$-structure, then a group homomorphism $Z\colon K(\DD)\to\C$ gives rise to a unique stability condition when the following two conditions are satisfied (\cite[Prop.\ 5.3]{Bridgeland:Stab}):
(i) $Z(\AA\setminus0)\subset\H$ ($Z$ is a \emph{stability function} on $\AA$);
(ii) Harder-Narasimhan filtrations exist for objects in $\AA$ with respect to $Z$.

Condition (i) means that, for all $0\neq A\in\AA$, the requirement $Z(A)\in\H$ gives a well-defined phase $\phi(A):=(1/\pi)\arg(Z(A))\in(0,1]$. This defines a notion of \emph{phase-stability} for objects in $\AA$, and so of (semi)stable objects of $\AA$. Then condition (ii) asks for the existence of finite filtrations for every object in $\AA$ in semistable ones with decreasing phases.

In particular, if $\AA$ is an abelian category of finite length (i.e., Artinian and Noetherian) with a finite number of
simple objects $\{S_0, \dots, S_m\}$, then any group homomorphism $Z\colon  K(\DD)\to\C$ 
with $Z(S_i)\in \H$ for all $i$ extends to a unique stability condition on $\DD$.
\end{Rem}

We give an improved criterion for the existence of Harder-Narasimhan filtrations:

\begin{Prop}\label{prop:HNFiltrationsDiscrete}
Let $\AA\subset\DD$ be the heart of a bounded $t$-structure on $\DD$ and let $Z\colon K(\DD)\to\C$ be a stability
function on $\AA$. Write $\PP'(1) \subset \AA$ for the full subcategory of objects with phase 1 with respect to
$Z$, and assume that:
\begin{itemize}
\item The image of $\Im(Z)$ is a discrete subgroup of $\R$.
\item For all $E\in\AA$, any sequence of subobjects
\[
0=A_0\subset A_1\subset\ldots\subset A_j\subset A_{j+1}\subset\ldots \subset E,
\]
with $A_j\in\PP'(1)$, stabilizes.
\end{itemize}
Then Harder-Narasimhan filtrations exist for objects in $\AA$ with respect to $Z$.
\end{Prop}
\begin{Prf}
We use the same ideas as in \cite[Prop.\ 7.1]{Bridgeland:K3}, and we want to apply \cite[Prop.\ 2.4]{Bridgeland:Stab}.

First of all notice that, if
\[
0\to A\to E\to B\to 0
\]
is an exact sequence in $\AA$, then
\[
0 \leq \Im Z(A)\leq \Im Z(E)\quad\text{and}\quad 0 \leq \Im Z(B)\leq \Im Z(E).
\]

Let
\[
\ldots\subset E_{j+1}\subset E_j\subset\ldots\subset E_1\subset E_0=E
\]
be an infinite sequence of subobjects of an object $E$ in $\AA$ with $\phi(E_{j+1})>\phi(E_j)$, for all $j$.
Since $\Im Z$ is discrete, there exists $N\in\NN$ such that
\[
0\leq \Im Z(E_n)=\Im Z(E_{n+1}),
\]
for all $n\geq N$.
Consider the exact sequence in $\AA$
\[
0\to E_{n+1}\to E_n\to F_{n+1}\to 0.
\]
Then, by additivity of $\Im Z$, we have $\Im Z(F_{n+1})=0$, for all $n\geq N$.
But this yields $\phi(F_{n+1})=1$, for all $n\geq N$ and so $\phi(E_{n+1})\leq\phi(E_n)$, a contradiction.
In this way, property $(a)$ of \cite[Prop.\ 2.4]{Bridgeland:Stab} is satisfied.

Let
\[
E=E_0\twoheadrightarrow E_1\twoheadrightarrow\ldots\twoheadrightarrow E_j\twoheadrightarrow E_{j+1}\twoheadrightarrow\ldots
\]
be an infinite sequence of quotients of $E$ in $\AA$ with $\phi(E_j)>\phi(E_{j+1})$, for all $j$.
As before, $\Im Z(E_n)=\Im Z(E_{n+1})$, for all $n\geq N$.
Consider the exact sequence in $\AA$
\[
0\to F_n\to E_N\to E_n\to 0,
\]
for $n\geq N$. Then $\Im Z(F_n)=0$, i.e., $F_n\in\PP'(1)$.
Hence we have an infinite sequence of subobjects of $E_N$ belonging to $\PP'(1)$, a contradiction.
Property $(b)$ of \cite[Prop.\ 2.4]{Bridgeland:Stab} is then verified and the proposition is proved.
\end{Prf}

A stability condition is called \emph{locally-finite} (see \cite[Sect.\ 5]{Bridgeland:Stab}) if there exists some $\epsilon>0$ such that, for all $\phi\in\R$, each quasi-abelian subcategory $\PP((\phi-\epsilon,\phi+\epsilon))$ is of finite length.
In this way $\PP(\phi)$ has finite length so that every object in $\PP(\phi)$ has a finite Jordan--H\"older filtration
into stable factors of the same phase. The set of stability conditions which are locally-finite will be denoted by
$\Stab(\DD)$. The stability conditions we consider also satisfy the additional conditions in the
definition given in \cite[Section 2]{Kontsevich-Soibelman:stability} (in particular the support
property, as discussed below). The local-finiteness condition will then be automatic.

The main result in \cite{Bridgeland:Stab} endows $\Stab(\DD)$ with a topology, induced by a metric $d(-,-)$ (see
\cite[Prop.\ 8.1]{Bridgeland:Stab} for the explicit form of $d$), in such a way it becomes a complex manifold whose
connected components are locally modeled on linear subspaces of $\Hom(K(\DD), \C)$ via the map $\ZZ$ sending a stability condition $(Z,\PP)$ to its central charge $Z$.

A connected component of $\Stab(\DD)$ is called \emph{full} if it has maximal dimension, i.e., it is modeled on the whole $\Hom(K(\DD), \C)$.
For a stability condition $\sigma=(Z,\PP)$ belonging to a full connected component (we will call $\sigma$ \emph{full}), we recall the statement of Bridgeland's deformation result.
In this case, the metric
\begin{equation*} 
\| W\|_{\sigma}:=\sup\left\{\frac{|W(E)|}{|Z(E)|}\colon E\text{ is }\sigma\text{-stable}\right\}
\end{equation*}
is finite, and thus defines a topology on $\Hom(K(\DD), \C)$.

\begin{Thm}[{\cite[Theorem 7.1]{Bridgeland:Stab}, \cite[Lemma 4.5]{Bridgeland:K3}}]
\label{thm:B-deform}
In the situation of the previous paragraph, let $0<\epsilon<1/8$.
Then, for any group homomorphism $W\colon K(\DD)\to\C$ with
\[
\|W-Z\|_{\sigma}<\sin(\pi\epsilon),
\]
there exists a unique (locally-finite) stability condition $\tau=(W,\QQ)\in\Stab^*(\DD)$ with $d(\sigma,\tau)<\epsilon$.
\end{Thm}
In particular this shows that the map $\ZZ \colon \Stab^*(\DD) \to \Hom(K(\DD), \C)$ is a local homeomorphism.

Let us also clarify the relation between full stability conditions in the situation of finite-rank
$K$-group and the support property introduced in \cite{Kontsevich-Soibelman:stability}.
More precisely, assume that $\DD$ is a triangulated category such that $K(\DD)/\mathrm{torsion}$ is a
finite-dimensional lattice, and choose a metric $\abs{\cdot}$ on $K(\DD)_\R$.
Then a stability condition $\sigma=(Z,\PP)\in\Stab(\DD)$ has the \emph{support property}\footnote{In the
notation of \cite[Section 1.2]{Kontsevich-Soibelman:stability}, we implicilty made the choices
$\Lambda = K(\DD)/\mathrm{torsion}$ and $\mathrm{cl}$ the projection.}
if there exists a constant $C>0$ such that
\[
C \abs{Z(E)} \ge \abs{E}
\]
for all $\sigma$-stable $E \in \DD$.

\begin{Prop}\label{prop:SupportProperty}
Assume that $K(\DD)$ has finite rank.
Then a Bridgeland stability condition $\sigma = (Z, \PP)$ is full if and only if it has the support property.
\end{Prop}
\begin{Prf}
Denote by $\abs{\cdot}^\vee$ the induced metric on $\Hom(K(\DD), \C)$.
A stability condition $\sigma=(Z,\PP)$ is full if and only if the semi-metric $\| \cdot \|_\sigma$ is finite.
Since $\Hom(K(\DD), \C)$ is finite-dimensional, this holds if and only if it is 
bounded by a multiple of $\abs{\cdot}^\vee$, i.e., if and only if there exists $C > 0$ such that, for any
$W \in \Hom(K(\DD), \C)$, we have
\[ 
\|W\|_\sigma \le C \cdot \abs{W}^\vee
\]
Hence if $\sigma$ satisfies the support property, then
\begin{equation*}
\begin{split}
\|W\|_\sigma &=\sup\left\{\frac{|W(E)|}{|Z(E)|}\colon E\text{ is }\sigma\text{-stable}\right\}\\
&\le C \cdot \sup\left\{\frac{|W(E)|}{\abs{E}}\colon E\text{ is }\sigma\text{-stable}\right\}\\
&\le C \cdot \abs{W}^\vee
\end{split}
\end{equation*}
and so $\sigma$ is full.

Conversely, assume that $\sigma$ does not satisfy the support property, i.e., there is a sequence
$E_n$ of $\sigma$-stable objects with $\abs{Z(E_n)} < \frac{\abs{E_n}}n$.
Let $W_n\in\Hom(K(\DD), \C)$ be such that $\abs{W_n}^\vee = 1$ and
$\abs{W_n(E_n)} = \abs{E_n}$. Then 
\[
\|W_n\|_{\sigma} \ge \frac{\abs{W_n(E_n)}}{\abs{Z(E_n)}}
> n \cdot \frac{\abs{E_n}}{\abs{E_n}} = n \abs{W_n}^\vee
\]
and so $\sigma$ is not full.
\end{Prf}

\begin{Rem}\label{rmk:GroupAction}
By \cite[Lemma 8.2]{Bridgeland:Stab}, we have a left action on $\Stab(\DD)$ by the autoequivalence group $\Aut(\DD)$, and a right action by $\widetilde{\GL}_2(\R)$, the universal cover of the matrices in $\GL_2(\R)$ with positive determinant. The first action is defined, for $\Phi\in\Aut(\DD)$, by $\Phi(Z,\PP)=(Z\circ\phi_*^{-1},\Phi(\PP))$, where $\phi_*$ is the automorphism induced by $\Phi$ at the level of Grothendieck groups.
The second one is the lift of the action of $\GL_2(\R)$ on $\Hom(K(\DD), \C)$ (by identifying $\C\cong\R^2$).
Notice, in particular, that the additive group $\C$ acts on $\Stab(\DD)$, via the embedding $\C\into\widetilde{\GL}_2(\R)$.
\end{Rem}

\section{Proof of the inequality for central charges}\label{app:ineq}

This appendix is a brief sketch of a complete proof of Observation (\ref{obs:inequality}) on
page \pageref{obs:inequality}: On the fundamental domain $\C \setminus \R_{\le 0}$, the 
functions $a(z), b(z)$ defined by equations \eqref{eq:ab-exp0} and \eqref{eq:ab-expinf}
satisfy the inequalities of Definition \ref{def:setG}.
The general idea is to deduce the inequalities from inequalities for real or imaginary parts
of holomorphic functions, which only need to be tested on the boundary of the fundamental domain. Note that the boundary,
expressed in $z$ and $\psi$ according to the convergence domains $\abs{z} \le \frac 1{27}$
and $\abs{\psi} < 1$ of our power series expansions,
consists of two copies of $z \in [-1/27, 0]$, with the two natural branch choices of $\ln z$,
and of the two ray segments
$\psi \in [0, 1] \cdot e^{2\pi i/3}$ and $\psi \in [0, 1] \cdot e^{4 \pi i/3}$. (We will refer to the boundary
segments by $\arg(z) = \pm \pi$ etc.)

\subsection*{Step 1.} We first show that $\Im(a)>0$.
For example by using the integral criteria, it can be shown easily that the series
$\sum_{n=1}^\infty \frac{(3n-1)!}{n!^3} \frac{1}{27^n}$
converges to a real number less than $\frac{1}{\pi\sqrt{3}}$ (the exact value is 0.129...).
Thus, for $\abs{z}<1/27$, $z\neq0$, $\abs{\arg(z)}<\pi$, we have
\[
\Im(a(z))\geq-\frac{1}{2\pi}\left(\ln\frac 1{27} +3 \sum_{n=1}^\infty \frac{(3n-1)!}{n!^3} \frac{1}{27^n}\right)>0.
\]
Passing to the $\psi$-coordinate, the inequality follows trivially on the two boundary rays
from \eqref{eq:ab-expinf} and the definition of $\varpi_1(\psi)$.

\subsection*{Step 2.} For $B=-\Im(b(z))/\Im(a(z))$ we have
$-1 <  B < 0$. More precisely, we will use $-1/2\leq B < 0$ for $\Im(z) \ge 0$ (i.e., for $2\pi/3 < \arg(\psi)\leq\pi$
or $0\leq\arg(z) < \pi$),
and $-1< B\leq -1/2$ for $\Im (z) \le 0$ (i.e., for $\pi\leq\arg(\psi)<4\pi/3$ or $-\pi<\arg(z)\leq0$).

To show this, first notice that $\Im(Z(\OO_{\P^2})(z))=\Im(b(z)) > 0$: indeed, for $\arg(z) = \pi$ or
$\arg(\psi) = \frac {2\pi}3$, it is trivial to check that $\Im(b(z)) = 0$. Similarly, the inequality holds strictly
for $\arg(z) = - \pi$ or $\arg(\psi) = \frac{4\pi}3$, and it also holds around $z = 0$. Thus the strict inequality
holds on the interior of the fundamental domain, and thus $B < 0$.

Similarly we can show $\Im(Z(\OO_{\P^2}(-1)[2])(z)) = \Im(b(z) - a(z)) < 0$ and thus $-1 < B$.
For the more precise statement, it is sufficient to look at the sign of
$\Im(Z(\Omega_{\P^2}(1)[1])(z)) = \Im(b(z) - \frac 12 a(z))$: once again the maximum principle shows that
$\Im(Z(\Omega_{\P^2}(1)[1])(z))\leq0$, for $2\pi/3 < \arg(\psi)\leq\pi$ and $0\leq\arg(z)< \pi$. This implies to $-\frac
12 \leq B$. The case $\Im (z) \le 0$ is analogous.

\subsection*{Step 3.} Finally, to check that the other two inequalities of Definition \ref{def:setG} are satisfied, we
show the following stronger statement:
\begin{equation}\label{eq:UpperSemiCircUtah}
-\Re(b(z))-B\Re(a(z))+\frac{B^2}{2}<\frac{3}{8}=\Delta_{-1/2},
\end{equation}
for all $-1/2\leq B\leq0$, in the region $\Im(z) \ge 0$ (and an analogous statement, which we will skip,
for $\Im(z) \le 0$). By the claims of the previous step, this will imply
Observation (\ref{obs:inequality}), as $\delta^{DP}_\infty(\mu) \ge \frac 12$, and as
$\Delta_B \ge \frac 38$ for $B \not\in \Z$. 
Also note that we only have to prove the inequality above for $B = 0$, and for $B = -\frac 12$.

First we note that $b(-1/27)=0$ (with the choice of $\arg(-1/27) = \pi$): this can be deduced from the monodromy.
From this, it follows that the series $\sum_{n=1}^\infty \abs{l_n} \frac{1}{27^n}$ (from the definition, it is easy to see that $l_n=(-1)^n\abs{l_n}$) converges to a real number less than 3 (the exact value is 1.558...).
From this we can deduce \eqref{eq:UpperSemiCircUtah} for the cases $\arg(z) = 0$ and $\arg(z)= \pi$:
Setting $B = - \frac 12$ and $B = 0$ yields
\begin{align*}
\tfrac 12 \Re(\omega_2(z))+\tfrac 18 &< \tfrac 38, \quad \text{and} \quad
\tfrac 12 \Re(\omega_2(z))-\tfrac 12 \Re(\omega_1(z))+\tfrac 14 <\tfrac 38.
\end{align*}
As $\Re(\omega_1(z))= \frac{1}{2\pi}\Im(\ln(z))$, this would follow from
\[
\tfrac 12 \Re(\omega_2(z))+\tfrac 18<\tfrac 38 \quad \text{if $z<0$, and } \quad
\tfrac 12 \Re(\omega_2(z))+\tfrac 14<\tfrac 38 \quad \text{if $z>0.$}
\]
Finally using the definition of $\omega_2(z)$, both inequalities become 
\begin{multline*}
-\frac{1}{8\pi^2} \left( (\ln\abs{z})^2 + 6 \ln\abs{z} \cdot \sum_{n=1}^\infty \frac{(3n-1)!}{n!^3} (-z)^n \right. 
			\left. + \sum_{n=1}^\infty \abs{l_n} (-z)^n \right) <\frac 18.
\end{multline*}
But the quantity on the left is at most
\[
-\frac{1}{8\pi^2} \left( \left(\ln\frac{1}{27}\right)^2 + 6 \ln\frac{1}{27} \cdot \sum_{n=1}^\infty \frac{(3n-1)!}{n!^3} \frac{1}{27^n}
			 - \sum_{n=1}^\infty \abs{l_n} \frac{1}{27^n} \right),
\]
which is smaller than $\frac 18$ by the estimate of the last term mentioned earlier.

For the cases $\arg(\psi) = \pi$ and $\arg(\psi) = \frac{2\pi}3$, we first observe that, for $0<\rho\leq1$ and
for any function $u\colon \Z_{>0}\to\{0,1\}$ with $u(1)=1$, we have
\begin{equation}\label{eqn:Estimatepsi}
\begin{split}
\frac 1{2\pi} &\sum_{n=1 \atop 3\nmid n}^\infty\frac{\Gamma\left(\frac n3 \right)}{\Gamma(n+1) \Gamma\left(1 - \frac n3 \right)^2} (3\rho)^n(-1)^{u(n)}\\
&\leq\rho\left(-\frac{3}{2\pi}\frac{\Gamma\left(\frac 13 \right)}{\Gamma\left(\frac 23 \right)^2}+\frac {1}{2\pi}\sum_{n=2 \atop 3\nmid n}^\infty\frac{\Gamma\left(\frac n3 \right)}{\Gamma(n+1) \Gamma\left(1 - \frac n3 \right)^2} 3^n\right)<0.
\end{split}
\end{equation}
Setting $B =0$ and $B = - \frac 12$, and using the definition of $\varpi_1(\psi)$, the needed inequalities
are
\[
\frac 38>\frac 13\Re(\varpi_2(-\rho))+\frac 13,
\]
for $\psi=-\rho$, and
\[
\frac 38>\begin{cases}
               \frac 16 \Re(\varpi_2(\rho e^{2\pi i/3}))+\frac 5{24}\\
	      \frac 23 \Re (\varpi_2(\rho e^{2\pi i/3})) + \frac 13
	      \end{cases},
\]
for $\psi=\rho e^{2\pi i/3}$, $0<\rho\leq1$.
But, by \eqref{eqn:Estimatepsi},
\[
\Re(\varpi_2(\rho e^{2\pi i/3})), \Re(\varpi_2(-\rho)) <0.
\]
Hence, \eqref{eq:UpperSemiCircUtah} holds also for $\arg(\psi) = \pi$ and $\arg(\psi) = \frac{2\pi}3$, and the proof of Observation \eqref{obs:inequality} is complete.

\bibliography{all}                      
\bibliographystyle{alphaspecial}     


\end{document}